\documentclass{gtmon_a}
\pdfoutput=1

\proceedingstitle{Heegaard splittings of 3--manifolds (Haifa 2005)}
\conferencestart{10 July 2005}
\conferenceend{19 July 2005}
\conferencename{Heegaard splittings of 3--manifolds}
\conferencelocation{Haifa}

\editor{Cameron Gordon}
\givenname{Cameron}
\surname{Gordon}

\editor{Yoav}
\givenname{Yoav}
\surname{Moriah}

\title{Heegaard splittings of knot exteriors}

\author{Yoav Moriah} 
\givenname{Yoav}
\surname{Moriah}
\address{Department of Mathematics\\
Technion\\\newline
Haifa 32000\\
Israel}
\email{ymoriah@tx.technion.ac.il}
\urladdr{}

\volumenumber{12}
\issuenumber{}
\publicationyear{2007}
\papernumber{07}
\startpage{191}
\endpage{232}

\doi{}
\MR{}
\Zbl{}

\arxivreference{math/0608137}



\received{4 August 2006}
\revised{30 March 2007}
\accepted{10 April 2007}
\published{3 December 2007}
\publishedonline{3 December 2007}
\proposed{}
\seconded{}
\corresponding{}
\version{}

\makeatletter
\def\cnewtheorem#1[#2]#3{\newtheorem{#1}{#3}[section]
\expandafter\let\csname c@#1\endcsname\c@theorem}

\let\xysavmatrix\xymatrix
\def\xymatrix{\disablesubscriptcorrection\xysavmatrix}
\AtBeginDocument{\let\bar\wbar\let\tilde\wtilde\let\hat\what}

\makeop{Aut}
\makeop{dim}
\makeop{gcd}
\makeop{genus}
\makeop{rank}
\makeop{min}
\makeop{max}

\theoremstyle{plain}
\newtheorem{theorem}{Theorem}[section]
\cnewtheorem{proposition}[theorem]{Proposition}
\cnewtheorem{lemma}[theorem]{Lemma}
\cnewtheorem{corollary}[theorem]{Corollary}
\cnewtheorem{claim}[theorem]{Claim}
\cnewtheorem{conjecture}[theorem]{Conjecture}
\theoremstyle{remark}
\cnewtheorem{question}[theorem]{Question}
\cnewtheorem{problem}[theorem]{Problem}
\cnewtheorem{definition}[theorem]{Definition}
\cnewtheorem{remark}[theorem]{Remark}
\cnewtheorem{example}[theorem]{Example}

\makeatother  

\def\reals{\mathbb{R}}
\def\naturals{\mathbb{N}}
\def\rationals{\mathbb{Q}}

\begin{document}

\begin{htmlabstract}
The goal of this paper is to offer a comprehensive exposition of the
current knowledge about Heegaard splittings of exteriors of knots in
the 3-sphere. The exposition is done with a historical perspective
as to how ideas developed and by whom. Several new notions are
introduced and some facts about them are proved.  In particular the
concept of a 1/n-primitive meridian. It is then proved that if a
knot K &sub; S<sup>3</sup> has a 1/n-primitive meridian; then nK = K#
&hellip; #K n-times has a Heegaard splitting of genus nt(K) + n
which has a 1-primitive meridian.  That is, nK is &micro;-primitive.
\end{htmlabstract}

\begin{abstract} 
The goal of this paper is to offer a comprehensive exposition of the
current knowledge about Heegaard splittings of exteriors of knots in
the $3$-sphere. The exposition is done with a historical perspective
as to how ideas developed and by whom. Several new notions are
introduced and some facts about them are proved.  In particular the
concept of a $1/n$-primitive meridian. It is then proved that if a
knot $K \subset S^3$ has a $1/n$-primitive meridian; then $nK = K\#
\cdots \#K$ $n$-times has a Heegaard splitting of genus $nt(K) + n$
which has a $1$-primitive meridian.  That is, $nK$ is $\mu$-primitive.
\end{abstract}

\begin{asciiabstract} 
The goal of this paper is to offer a comprehensive exposition of the
current knowledge about Heegaard splittings of exteriors of knots in
the 3-sphere. The exposition is done with a historical perspective as
to how ideas developed and by whom. Several new notions are introduced
and some facts about them are proved.  In particular the concept of a
1/n-primitive meridian. It is then proved that if a knot K in S^3 has
a 1/n-primitive meridian; then nK = K#...#K, n-times has a Heegaard
splitting of genus nt(K) + n which has a 1-primitive meridian.
That is, nK is mu-primitive.
\end{asciiabstract}

\maketitle

\section {Introduction}

The goal of this survey paper is to sum up known results about
Heegaard splittings of knot exteriors in $S^3$ and present them with
some historical perspective. Until the mid 80's  Heegaard splittings
of $3$--manifolds and in particular of knot exteriors were not well
understood at all. Most of the interest in studying knot spaces,
up until then, was directed at various knot  invariants which had a
distinct algebraic flavor to them. Then in 1985 the remarkable work of
Vaughan Jones turned the area around and began the era of the modern
knot invariants \`{a} la the Jones polynomial and its descendants and
derivatives. However at the same time there were major developments
in Heegaard theory of $3$--manifolds in general and knot exteriors
in particular. For example the invention of the notions of strongly
irreducible and weakly reducible Heegaard splittings by Casson and
Gordon in \cite{CG} and various techniques to deal with when two Heegaard
splittings are the same by Boileau and Otal in \cite{BO}. Further work
was done on the subject of distinguishing Heegaard splittings by Lustig
and the author.  A multitude of results were obtained, at that time,
by many mathematicians, such as the Japanese school led by Kobayashi,
Morimoto and Sakuma. Other results were obtained by  Scharlemann,
Schultens, Thompson and many others.\footnote{ My apologies to those I
have omitted in this short list. I hope to rectify this in the body of
the paper.}  It is these results that will be surveyed and discussed.

Many of the definitions  and the results which are brought here have been
defined and and applied to general $3$--manifolds. However they will be
treated here in the more narrow context of $3$--manifolds $E(K) = S^3 -
N(K)$ which are {\it exteriors of knots in $S^3$}.

When studying Heegaard splittings of $3$--manifolds there are  some basic
topics  that we study:

\begin{description}
\item[Genus] Determining the minimal genus of $E(K)$.
\item[Classification] This topic splits into the following:
\begin{enumerate}
\item Finding  inequivalent Heegaard splittings.
\item Determining all Heegaard splittings up to equivalence (see \fullref{preliminaries}).
\item Deciding when Heegaard splittings are strongly irreducible or
weakly reducible
\item Understanding non-minimal genus Heegaard splittings. 
\end{enumerate}
\item[Additivity] Let $g(M)$ denotes the genus of the manifold
$M$. Suppose $K_1,K_2 \subset S^3$  are knots and we given $g_1 =
g(E(K_1))$ and $g_2 = g(E(K_2))$.  What can we say about $g(E(K_1\#K_2))$? 
\end{description}

For definitions, terminology and preliminary facts see \fullref{preliminaries}.

In general the methods that are used to study the above topics are
basically topological and algebraic.  The topological methods make use
of the quiver of the now standard techniques that were developed in
the late eighties and early nineties, for example, the ideas of  {\it
thin position} introduced by Gabai  to get information about Heegaard
splittings of knot spaces and the application of {\it  Cerf Theory} by
Rubinstein and Scharlemann in order to compare two different Heegaard
splittings, also the notions of {\it strongly irreducible,  weakly
reducible} and the {\it Rectangle Condition} introduced by Casson and
Gordon. These are the technical tools that are responsible for many of
the topological results.  Unfortunately it is  beyond the scope of this
paper to describe these techniques in detail.

There are  also some algebraic ideas developed by M. Lustig and the
author that enable us to determine the rank of groups in many cases,
in this case the fundamental group of $E(K)$. These techniques also
distinguish generating systems up to {\it Nielsen Equivalence} and this
was used to distinguish Heegaard splittings up to isotopy.

One of the problems that arose writing this survey was placing  the many
results in the different sections and sub-sections. Many of the results
play a role in more than one  aspect of understanding the Heegaard
splittings of knot exteriors so they are mentioned not necessarily in
the obvious order.

\subsection*{Acknowledgments}
I wish to thank the referee for numerous remarks and a meticulous job.
This work was partially supported by grant 2002039 from the US--Israel
Binational Science Foundation (BSF), Jerusalem, Israel.

\section{Preliminaries}\label{preliminaries}

A {\it compression body} is a $3$--manifold $V$ obtained from a surface $S$
cross an interval $[0,1]$  by attaching a finite number of $2$--handles
and $3$--handles to $S\times\{0\}$.  The component $S\times\{1\}$ of  the
boundary will be denoted by $\partial_{+}V$, and $\partial V \setminus
\partial_{+}V$ will be denoted by  $\partial_{-}V$.  The trivial cases
where $V$ is a handlebody or $V = S \times [0,1]$ are allowed.

Let $K\subset S^3$ be a knot. A {\it Heegaard splitting } for a knot
exterior $E(K)$ is a decomposition  $E(K) = V \cup_S W$, where $V$ is
a compression body with $\partial_+V = S$ and $\partial_-V = \partial
N(K)$. Two Heegaard splittings for a given $E(K)$ will be called
equivalent up to  homeomorphism (isotopy) if there is a homeomorphism
(isotopy) $h\co E(K) \rightarrow E(K)$ such that $h(S) = S$. A regular
neighborhood will be denoted by $N(~)$.

Historically when studying knots people used the notion of unknotting
tunnels which is essentially equivalent  to that of  genus:

\begin{definition}
Given a knot $K \subset S^3$ a collection of disjoint arcs $t_1, \dots,
t_n$ properly embedded in $E(K)$ will be called an  {\it unknotting tunnel
system} if $E(K) - \cup \{ N(t_i)\}$ is a  handlebody.  Note that $N(K)
\cup (\cup \{ N(t_i)\})$ is always a handlebody and hence a tunnel system
determines a Heegaard splitting of $S^3$ with the property that $K$
is contained as a core curve of one of the handlebodies..  The minimal
cardinality of any such tunnel system for the given knot $K$ is called
the {\it tunnel number} of $K$ and is denoted by $t(K)$.
\end{definition}

It follows immediately that given a tunnel system with $n$ tunnels
the genus, of the Heegaard splitting of $E(K)$ it determines, is $g =
n +1$. In particular a minimal tunnel system determines a minimal genus
Heegaard splitting of $E(K)$ and $g(E(K)) = t(K) + 1$.

As shall be seen it is convenient for many purposes  to consider knots
in $S^3$ as $2n$--plats. For convenience the definition is stated below.

\begin{definition}
\label{2nplat}
Given a regular projection of a braid on $2n$ strings we can cap off
consecutive pairs of strings on the top and on the bottom of the braid
by small arcs, called {\it bridges} to get a projection of a knot or
link in $S^3$. Such a projection is called a {\it $2n$--plat} (see Burde
and Zieschang \cite{BZ1}).
\end{definition}

\begin{remark}
Note that:
\begin{enumerate}
\item All knots and links  $ K \subset S^3$ have such projections. This
follows from the theorem due to Alexander \cite{Al} that all knots and
links have a representation as closed braids. There is an algorithm to
obtain a braid presentation from a knot projection.

\item  If a knot $K \subset S^3$ has a $2n$--plat projection then the
bridge number $b(K)$ of $K$ is bounded above by $n$.

\item A $2n$--plat projection determines two sets of unknotting tunnel
systems for $K$ by connecting consecutive bridges by $n{-}1$ small arcs
(see Lustig and Moriah \cite{LM2}).

\item A $2n$--plat projection determines two Heegaard splittings of genus
$n$ for $E(K)$.

\item A $2n$--plat projection determines two sets of $n$ generators for
$\pi_1(E(K))$:

Let ${\hat x} = \{x_1, \ldots, x_n\}$ be a collection of
small circles around the top bridges and ${\hat y} = \{y_1, \ldots, y_n\}$
be a collection of small circles around the bottom bridges. Connect both
${\hat x}$ and $\hat y$ to a chosen base point. Then the curves $\hat x$
and $\hat y$ are representatives of generators.

\item Such a projection determines two presentations for $\pi_1(E(K))$:

Denote the top tunnels  by $\tau_1, \ldots, \tau_{n - 1}$ and the
bottom tunnels by $\eta_1, \ldots, \eta_{n - 1}$.  Consider a small
regular neighborhood of the tunnels to obtain two collections of $n
- 1$ $1$--handles. The cocore disks of the $\eta_1, \ldots, \eta_{n -1}$
$1$--handles determine a set of relations for the top generators and the
cocore disks of the $\tau_1, \ldots, \tau_{n - 1}$ $1$--handles determine
a set of relations for the bottom generators \cite{LM2}.
\end{enumerate}
\end{remark}

One other feature of Heegaard splittings which is relevant to this
discussion is the notions of strongly irreducible and weakly reducible
which are due to Casson--Gordon \cite{CG} and are defined as follows:

\begin{definition} \label{reducibility}
A Heegaard splitting $(V,W)$ of an irreducible  $3$--manifold $M$
will be called {\it weakly reducible} If there are essential disks $D
\subset V$ and $E \subset W$ such that  $\partial D \cap \partial
E = \emptyset$. Otherwise $(V,W)$ will be called {\it strongly
irreducible}. If there are disks such that $\partial D \cup \partial E $
is a single point $p$ we say that  $(V,W)$ is {\it stabilized}.
\end{definition}

The following is a well known definition due to Harvey which we bring
for completeness:

\begin{definition}\label{curvecomplex}
Given a surface $S$ of genus $g \geq 2$ let  ${\mathcal C}_S$ denote
the {\it curve complex}  of $S$ defined as follows:
\begin{enumerate}
\item The  vertex set,  ${\mathcal V}_S$ of $\mathcal C_S$ is the set
$\{[\gamma]\}$ of isotopy classes of essential simple closed curves
$\gamma$ on $S$.

\item A set of $n$  vertices $\{v_1, \dots, v_n\}$ represented by a set
of $n$ disjoint curves will define an $(n{-}1)$--simplex.
\end{enumerate}
\end{definition}

Note that $\mathcal C_S$ is  connected, $\dim(\mathcal C_S) =
3g-4$, where $g$ is the genus of $S$ and that it is not locally
finite. On the $1$--skeleton of $\mathcal C_S$ there is a metric
$d_{\mathcal{C}(S)}(\cdotp,\cdotp)$ defined by  setting the length of
every edge to be one.

\begin{definition}
\label{diskcomplexes}
Given a $3$--manifold $M$ with a Heegaard splitting $(V,W)$, where $V$ and
$W$ are  compression bodies and a  Heegaard surface $S $, the collection
of vertices which correspond to essential curves in $S$ which bound
disks in $V$ and $W$ respectively, define sub-complexes  $\mathcal D_V$
and $\mathcal D_W$ of $\mathcal C_S$.  They are called the {\it disk
complexes} of the corresponding compression bodies.
\end{definition}

Note that the complexes $\mathcal D_V$ and $\mathcal D_W$ of $\mathcal
C_S$ are connected.  The following definition is due to Hempel
~\cite{He1}:

\begin{definition} \label{distance}
Given a Heegaard splitting $(V, W)$ for a $3$--manifold $M$ with $
\partial_{+} V = S =  \partial_{+} W$, we define the {\it distance of
the Heegaard splitting}, denoted by $d(V, W)$,  as follows: $$d(V,W)
= \min \{d_{\mathcal{C}(S)}(\partial D, \partial E) | D \subset V,
E \subset W\} $$ where $D \subset V, E \subset W$ are essential disks
and $d_{\mathcal{C}(S)}( \cdotp,\cdotp )$ is the distance in  $\mathcal
C_S$. In other words, it is the distance between  $\mathcal D_V$ and
$\mathcal D_W$ in $\mathcal C_S$.
\end{definition}

Note that for a given Heegaard splitting $(V,W)$, $d(V,W) = 0$ is
equivalent to $(V,W)$ being reducible, $d(V,W) = 1$ is equivalent to
$(V,W)$ being weakly reducible but irreducible and $d(V,W) \geq 2 $
is equivalent to $(V,W)$ being strongly irreducible.

\begin{remark}
Recall that all knot spaces are irreducible $3$--manifolds. Hence Heegaard
splittings  of distance zero are stabilized.
\end{remark}

For standard terminology and definitions see Burde and Zieschang
\cite{BZ1}, Rolfsen \cite{Ro}, Hempel \cite{He}  and  Jaco \cite{Ja}. Note
also that a lot of the results apply to Heegaard splittings of
$3$--manifolds which are not knot spaces.

\section{Determining the genus of $E(K)$}
\label{genus}

In theory, the problem of determining the genus for a knot space $E(K)$
is ``solved" by the Rigidity Theorem proved by K Johannson in \cite{Jo1}:

\begin{theorem}[Johannson ~\cite{Jo1}]
Let $M$ be a Haken $3$--manifold with or without boundary and  without a
non-trivial, essential Stallings fibration. Then the set of all isotopy
classes of Heegaard surfaces in $M$ of any given genus is finite and
constructable. In particular the genus of $M$ can be determined.
\end{theorem}

Since knot spaces are Haken then for non-fibered knots we are done by the
above theorem.  However, in practice, for a random knot the algorithm
suggested is quite impractical.  Upper bounds which are not sharp are
readily computable. When a knot $ K \subset S^3$ is given as a knot
projection into a  plane $P$ with say, $c_P(K)$ crossings, the number
$c_P(K)$ is  an upper  bound on its tunnel number. As, if we insert a
tunnel connecting the two strings involved in each crossing we obtain
a graph whose complement  in $S^3$ is  isotopic to the  complement of a
planar graph and hence is a handlebody.  However computing the crossing
number $C(K)$ of an arbitrary knot is in general very difficult. For
example, computing the unknotting number for all knots with 10 or fewer
crossings has only recently been accomplished by Gordon and Luecke
\cite{GL} and has required work of Ozsvath--Szabo  \cite{OS}.

Since obtaining upper bounds which are clearly not sharp is quite easy
the real question is, can we improve the upper bound? This question is
on an entirely different level of difficulty for general knots.

Looking for a lower bound we note that a Heegaard splitting $(V,W)$
of $E(K)$ of genus $g$ induces a presentation :
$$\pi_1(E(K)) = \langle x_1, \ldots, x_g | R_1, \ldots, R_{g - 1}\rangle$$
where $ \{x_1, \ldots, x_g\}$
represent generators of $\pi_1(W)$ given by an appropriate choice of  a
spine for the handlebody $W$. Thus the rank of $\pi_1(E(K))$ is clearly
a lower bound for the genus of any Heegaard splitting.

\begin{remark}\label{rank+ genus}
It is still unknown if, for knots in $S^3$, the rank of $\pi_1(E(K))$
is equal to $g(E(K))$. So we do not know if the determination
of $\rank(\pi_1(E(K)))$ will give a sharp lower bound for the
genus. In general there are examples by Boileau and Zieschang \cite{BZ} of Seifert
fibered spaces $M$ for which $\rank(\pi_1(M)){=}2$ but  $g(M){=}3$. There are results by Schultens--Weidmann \cite{SW}
for graph manifolds $M$ with $\rank(\pi_1(M)){=}3$ and $g(M){=}4$. The
question of whether rank equals genus for  hyperbolic $3$--manifolds is
an important open question. There are some  recent partial  results of
I Agol (unpublished), D Bachmann, D Cooper and M White \cite{BCW}
and last but not least J Souto \cite{So}.
\end{remark}

For hyperbolic knots there are results regarding the rank of the
fundamental group of the knot space. In \cite{KW} I Kapovitch and R
Weidmann prove:

\begin{theorem}[Kapovitch--Weidmann~\cite{KW}]
If $K\subset S^3$ is a hyperbolic knot then there exists an algorithm
which computes the rank of $\pi_1(E(K))$.
\end{theorem}

However, quoting Kapovitch and Weidmann  themselves: ``we should stress
that the (above) theorem is an abstract computability result. The nature
of the proof is such that it cannot provide any complexity bound on the
running of the algorithm."

The moral from the above theorems is that a general algorithm to
compute either the genus of the Heegaard splitting or the rank of the
fundamental group will be very hard to implement.  Hence we need to
yield in generality in order to gain computability.

Below we present an algebraic method (devised in Lustig--Moriah
\cite{LM2}) for computing the rank of a {\it general} finitely generated
group. It can be {\it readily} implemented to computing the rank of the
fundamental group of a very large class of knots and links and the only
condition that is required is that the knot/link be given in a certain
$2n$--plat projection. If a combinatorial condition, which is decided
by counting, is satisfied by the $2n$--plat, then the rank of the group
is determined. The issue then becomes:

\begin{question}
Does a knot $K \subset S^3$ have ``such" a $2n$--plat projection? (See 
Lustig and Moriah \cite[Theorems~0.1 and~0.4]{LM3}.)
\end{question}

A presentation, as above, determines a resolution $\varphi\co F(X_1, \ldots,
X_g) \rightarrow \pi_1(E(K))$ given by the homomorphism $\varphi\co X_i
\rightarrow x_i , i = 1, \ldots, g$, where $F_g = F(X_1, \ldots, X_g)$ is the
free group on $\{X_1, \ldots, X_g\}$.  Let $\mathbb{Z} G$ denote the group
ring of a group $G$ and  $\smash{\frac{\partial}{\partial X_i}}\co \mathbb{Z}
F(X_1, \ldots, X_g) \rightarrow \mathbb{Z} F(X_1, \ldots, X_g)$  be the $i$th
Fox derivative with respect to the generating system $\hat X = \{X_1, \ldots,
X_g\}$ (for definitions and discussion of the {\it Fox derivative} see
Fox \cite{Fo} and Lyndon--Schupp \cite{LS}). The relations $R_1, \ldots, R_{g - 1} $ are words
in $F(X_1, \ldots, X_g)$ so we can take their Fox derivatives. Furthermore
we can extend the homomorphism $\varphi$ to a ring homomorphism $\varphi
\co  \mathbb{Z} F(X_1, \ldots, X_g) \rightarrow \mathbb{Z} \pi_1(E(K))$.
Now let $I_{\hat X}$ denote the two sided ideal in $\mathbb{Z} \pi_1(E(K))$
generated by $\bigl\{\varphi \bigl(\smash{\frac{\partial R_j}{\partial
X_i}}\bigr)\bigr\}$ for $i = 1,
\ldots, g$ and $j = 1, \ldots, g - 1$, and let $\mathbb{Z} \pi_1(E(K))/ I_{\hat X}$
denote the quotient ring.  Let $R$ be some commutative ring with a unit,
$m \in  \mathbb{N}$ and $\mathbb{M}_m(R)$ be the $m{\times}m$ matrix
ring with entries in $R$.

We can now state the following theorem which is an immediate consequence
of Lustig--Moriah \cite[Lemma~2.1]{LM2}:

\begin{theorem}[Lustig--Moriah \cite{LM2}] \label{rank}
Suppose that for some $m \in \mathbb{N}$ and some commutative
ring with a unit $R$, there is non-trivial representation $\sigma\co 
\mathbb{Z} \pi_1(E(K))/ I_{\hat X} \rightarrow \mathbb{M}_m(R)$.
Then $\rank(\pi_1(E(K))) = g$.
\end{theorem}

Though at first glance implementing this theorem seems impossible due
to the impenetrable nature of  $\mathbb{Z} \pi_1(E(K))/ I_{\hat X}$,
it turns out that in practice  it is quite easy to apply. For example,
for knots with a presentation for the fundamental group of $E(K)$ coming
from  a $2n$--plat projection.  By setting  $m = 1$ and $R = \mathbb{Z}
[t, t^{-1}]$ there  are quite a few results which were computed by hand
and which will be described below.

Note that for knots $K \subset S^3$ in a $2n$--plat projection we have
the following inequalities:
$$\rank(\pi_1(E(K)))\leq g(E(K))  \leq b(K)  \leq n$$
Here $b(K)$ denotes the classical bridge number of $K$ with respect to
a $2$--sphere.

\fullref{rank} can be applied to a large class of knots and links in
$S^3$ defined in Lustig and Moriah \cite{LM3} which are called {\it generalized Montesinos
knots/links}. These include in particular all $2$--bridge knot/links,
all Montesinos knots (for a definition see Burde and Zieschang
\cite[pages 196--207]{BZ1})
and many other much more complicated knots/links as in Lustig and Moriah
\cite[Theorem~5.6]{LM3}. A generalized Montesinos knot/link is a combination via
braids of a collection of $2$--bridge  knot/links each with invariants
$\alpha_{i,j}/\beta_{i,j}$. If we set $\alpha = \gcd\{\alpha_{i,j}\}$
then we have:

\begin{theorem}[Lustig--Moriah {{\cite[Theorems~0.1(1) and~0.4(1)]{LM3}}}]\quad
\label{montesinos}
\begin{enumerate}
\item If $\alpha \neq 1$  then $\rank(\pi_1(E(K))) =  g(E(K))  =  b(K)  = n$.

\item If $K_1$ and $K_2$ are two generalized Montesinos knots and $K =
K_1\# K_2$ and if $\alpha \neq 1$  then $t(K) = t(K_1) + t(K_2)$.
\end{enumerate}
\end{theorem}

\begin{remark}
Let $K$ be a generalized Montesinos link with say $d$ components.
The above results on the genus of the knot/link space $E(K)$ can be
extended to closed manifolds obtained by $p_i/q_i , i = 1, \ldots, d $,
with $p_i$ even, surgery on $K$. See also \cite[Theorem~5.6]{LM3}.
\end{remark}

\begin{remark}
The above Theorems \ref{rank}, \ref{montesinos} and
\cite[Theorem~5.6]{LM3} compute the \emph{minimal} genus of $E(K)$ when they can
be applied. This is not a coincidence.  It follows from Lustig and Moriah
\cite[Theorem~2.1]{LM2} that there is no non-trivial representation
$\sigma\co 
\mathbb{Z} \pi_1(E(K))/ I_{\hat X} \rightarrow \mathbb{M}_m(R)$ if the
number of generators in the presentation is not minimal.
\end{remark}

\section{Inequivalent Heegaard splittings of knot spaces}
\label{inequivalent}

When trying to classify the Heegaard splittings of the exterior $E(K)$
of a given (say by a projection)  knot $K \subset S^3$, two problems arise
immediately. The first is how to obtain all the Heegaard splittings. The
second is how to show that they are different if they happen to be of
the same genus and one does not know that they are different by obvious
reasons. For example, one is strongly irreducible and the other weakly
reducible. (The issue of  strongly irreducible and  weakly reducible will
be addressed in \fullref{Strong+weak}.)  Heegaard splittings for $E(K)$
are equivalent to unknotting tunnel systems so the first problem can be
approached by trying to find inequivalent systems of unknotting tunnels.
The main tool for the second problem is the notion of {\it Nielsen
equivalence}. We will discuss this first. We define:

\begin{definition}\label{NE}
Let $\bar x = \{x_1, \ldots, x_n\}$ and $\bar y = \{y_1, \ldots, y_n\}$
be two systems of generators for a group $G$.  Denote by $F(X)$ and
$F(Y)$ the free groups on bases $X = \{X_1, \ldots, X_n\}$ and $Y = \{Y_1,
\ldots, Y_n\}$ respectively. Let $\beta_x$ and $\beta_y$ be the canonical
epimorphisms $F(X) \rightarrow G$ given by $X_i \rightarrow x_i$ and
$F(Y) \rightarrow G$ given by $Y_i \rightarrow y_i$. We will say that the
generating systems $\bar x$ and $\bar y$ are {\it Nielsen equivalent}
if there is an isomorphism $\alpha\co F(Y) \rightarrow F(X)$ such that $\beta_y =   \beta_x \circ \alpha$.
\end{definition}

\begin{remark}\label{freegp}
It is clear from the definition that all generating systems of cardinality
$n$ of a free group of rank $n$ are Nielsen equivalent.
\end{remark}

Suppose now that $E(K)$ has two genus $g$ Heegaard splittings $(V,W)$ and
$(P,Q)$.  There are spines of $W$ and $Q$ determining sets of generators
$\bar x = \{ x_1, \ldots, x_g\}$ and $\bar y = \{ y_1, \ldots, y_g\}$  for
$\pi_1(E(K))$. Assume that $(V,W)$ and $(P,Q)$ are isotopic. Then the
isotopy takes $Q$ to $W$, say. Hence the spine of $Q$ which determines a
set $\bar y = \{y_1, \ldots, y_g\}$  will be taken to a spine of $W$ which
determines a different set $\bar x' = \{x'_1, \ldots, x'_g\}$ of generators
for $\pi_1(W)$. Since $\pi_1(W)$ is a free group the sets $\bar x$
and $\bar x'$ are Nielsen equivalent. Hence the generators $\bar x$
and $\bar y$ are Nielsen equivalent.

If the Heegaard splittings are equivalent by a homeomorphism then the two
generating sets are Nielsen equivalent up to an automorphism of $G$. (See
for example Lustig and Moriah \cite{LM1,LM3,LM2}.)  It follows from
the above discussion that if one can show that the generating sets of
$\pi_1(E(K))$  which are determined are Nielsen inequivalent then the
corresponding Heegaard splittings  are not isotopic.
 
\subsection{Heegaard splittings of tunnel number one knots} \label{t(K)=1}

Tunnel number one knots/links have genus two Heegaard splittings and
hence the fundamental groups of their exteriors are generated by two
elements. For these groups we have the following theorem due to Nielsen
\cite{Ni}:

\begin{theorem}[Nielsen~\cite{Ni}]
Let $F$ be the free group with basis $\{x,y\}$ and let $w$ be the cyclic
word determind by the commutator $[x,y]$ . Then every automorphism of $F$
carries $w$ to $w$ or $w^{-1}$.
\end{theorem}

As an immediate consequence we have:

\begin{theorem} \label{commutator}
Let $G$ be a two generator group. Then the commutator of the generators
is, up to conjugacy and inverses, an invariant of the Nielsen class of
the generating set.
\end{theorem}

\begin{remark}
Any set of words in a free group $F$ of rank greater than two which is
invariant under $\Aut F$ is infinite. So there is no hope for extending
the above method to groups with rank bigger than two (see Lyndon and
Schupp \cite[Proposition~5.2]{LS}).
\end{remark}

\subsubsection{Torus knots} \label{torusknots}

The above theorem tells us that  if we can show that commutators or
their inverses of the generators in two generating sets determined by
two genus two Heegaard splittings are not conjugate then the Heegaard
splittings are not isotopic. This idea was very successful when studying
torus knots. The fundamental group of the exterior of a torus knot is
in general a $\mathbb{Z}$ extension of a triangle group which embeds in
$SL_2(\reals)$. One can find cores for the handlebodies which correspond
to elliptic elements in the image of the group in $SL_2(\reals)$. The
commutators are hyperbolic elements and hence have a translation length
which is invariant under conjugation and inverses.  Thus if one can show
that the translation lengths of commutators corresponding to different
generating sets are different the corresponding Heegaard splittings are
non-isotopic. This was used to obtain the following:

Let $p, q \in \naturals$ such that $\gcd(p,q) = 1$ and  $0 < p < q$.
Exteriors of torus knots $T(p,q)$ are Seifert fibered spaces over the
disk $D$ with two exceptional fibers $u, v$. They are denoted by $D(-r/p,
s/q)$ where $r,s \in \mathbb{Z}$ such that $ps - rq = 1$ and the Seifert
invariants of $u$ and $v$ are  $-r/p$ and $ s/q$ respectively. These
knots have three unknotting tunnels $\tau_u, \tau_v$ and $\tau_m$. The
tunnel $\tau_u$ (resp. $\tau_v$) is obtained by connecting the fiber $u$
(resp. $v$) by a small arc on the base space to the boundary.  The tunnel
$\tau_m$ is a small arc connecting $u$ to $v$ on $D$. A Heegaard surface
for such a Seifert  fibered space will be called  {\it vertical} if it
is the boundary of a regular neighborhood of $u \cup \tau_m \cup v$  or
a regular neighborhood of $\partial E(T(p,q)) \cup \tau_u$ or $\partial
E(T(p,q)) \cup \tau_v$. This definition is extended to Heegaard splittings
of general Seifert fibered spaces in Moriah and Schultens \cite{MS1}.

It was shown by Boileau--Otal in \cite{BO}, and by Boileau--Rost--Zieschang
in \cite{BRZ}, that all genus two Heegaard splittings of these
spaces are isotopic to the {\it vertical} ones. The method given by
\fullref{commutator} was used  in Moriah \cite{M,Mr0} to distinguish
the vertical Heegaard splittings and to obtain:

\begin{theorem}[Moriah~\cite{M,Mr0}]
For a torus knot $T(p,q)$ we have:
\begin{enumerate}
\item If $p \equiv \pm1$ mod $q$ then $\tau_m$ is isotopic to $\tau_u$
and $E(T(p,q))$ has two inequivalent vertical Heegaard splittings.

\item If $q \equiv \pm1$ mod $p$ then $\tau_m$ is isotopic to $\tau_v$
and $E(T(p,q))$ has two inequivalent vertical Heegaard splittings.

\item Otherwise $E(T(p,q))$ has exactly three inequivalent genus two
Heegaard splittings and they are vertical.
\end{enumerate}
\end{theorem}

The Heegaard splittings of torus knot exteriors of genus greater than
two were dealt with in Moriah and Schultens \cite{MS1}.  They are
stabilizations of the vertical minimal genus two splittings described
above.

\subsubsection{$2$--bridge knots} \label{$2$--bridge knots}

The study of $2$--bridge knots has a long history. For an excellent
exposition see Burde and Zieschang~\cite{BZ1}. They were  first
classified by Schubert in \cite{Sh} who showed that they are classified by
a rational number $\alpha/ \beta, -\alpha < \beta < \alpha $ and two knots
corresponding to $\alpha/ \beta$ and $\alpha'/ \beta'$ are equivalent if
and only if $\alpha = \alpha'$ and either $\beta = \beta'$ or $\beta\beta'
\equiv $1 mod $\alpha$. To a continued fraction expansion of the number
$\alpha/ \beta$ corresponds a rational tangle so $2$--bridge knots and
links can be described  by a $4$--plat  and the top and bottom tunnels of
the $4$--plat give immediate candidates for genus two Heegaard splittings.

The Heegaard  splittings determined by the top and bottom unknotting
tunnels for $2$--bridge knots  which are non-amphicheiral and
non-palindromic were shown to be inequivalent by Funcke \cite{Fu}. He
showed by a direct argument that the generating systems determined by
these Heegaard splittings are Nielsen inequivalent.

This was improved upon by  Bleiler--Moriah \cite{BM} for all $2$--bridge
knots and the distinction is up to homeomorphism. The method used was
geometric, making use of the fact that tunnel number one knots are
strongly invertible. Showing that the involutions are not equivalent
implies that the tunnels and hence the Heegaard splittings are
inequivalent as well.

\begin{remark}\label{AlgnotGeom}
In particular it is shown in \cite[Theorem~5.5]{BM} that for amphicheiral
$2$--bridge knots the Heegaard splittings determined by the top and bottom
unknotting tunnels are inequivalent though the generating systems for
$\pi_1(E(K))$ are equivalent.
\end{remark}

There are altogether six ``more or less" obvious unknotting tunnels. These
are the {\it top } and {\it bottom} tunnels $\tau_1, \tau_2$ respectively,
plus two arcs $\rho_1, \rho'_1$ with end points on the left (right)
top bridge and linking the right (left) top bridge once. Also two
tunnels $\rho_2, \rho'_2$ with the same construction for the bottom
bridges. These unknotting tunnels were classified by Kobayashi in
\cite{Ko1} up to homeomorphism and by Morimoto--Sakuma in \cite{MS3}
up to isotopy as follows:

Note that unknotting tunnels are equivalent if and only if the Heegaard
splittings are equivalent. The symbol $\cong$ will denote isotopy.

\begin{theorem}[Morimoto--Sakuma~\cite{MS3}]
Let $K(\alpha/\beta) \subset S^3, \alpha/\beta\in \rationals$, $\alpha$
odd, be a  $2$--bridge knot.  We can choose $\beta$ even. Then
\begin{enumerate}

\item $\tau_1 \cong \tau _ 2$  and $\rho_2, \rho'_2 \cong \rho _1 $
and $\rho_2, \rho'_2 \cong \rho'_ 1$ if and only if $\beta \equiv \pm1$
mod $\alpha$.

\item $\rho_1, \rho'_1 \cong \tau _ 1$  and $\rho_2, \rho'_2 \cong \tau _
2$ if and only if $\alpha = 3$.

\item $\rho_2, \rho'_2 \cong \tau _ 1$  if and only if $\beta \equiv \pm2$
mod $\alpha$.

\item $\rho_1, \rho'_1 \cong \tau _ 2$  if and only if $\beta \equiv
\pm2^{-1}$ mod $\alpha$.

\item $ \rho'_1 \cong \rho_ 1$  if and only if $\beta \equiv \pm1$
mod $\alpha$ and $\beta \equiv \pm2^{-1}$ mod $\alpha$.

\item $ \rho'_2 \cong \rho_ 2$  if and only if $\beta \equiv \pm1$
mod $\alpha$ and $\beta \equiv \pm2$ mod $\alpha$.
\end{enumerate}
\end{theorem}

Then Kobayashi,  in \cite{Ko2}, using the idea of a {\it labeled graphic}
given by Cerf theory and developed by Rubinstein and Scharlemann \cite{RS}
as a tool to compare  Heegaard splittings showed that:

\begin{theorem}[Kobayashi~\cite{Ko2}]
Every unknotting tunnel for a non-trivial $2$--bridge knot is isotopic
to one of the above six tunnels.
\end{theorem}

In a later beautiful paper Kobayashi,  \cite{Ko4} again using the method
of \cite{RS}, was able to classify all the non-minimal genus Heegaard
splittings of exteriors of $2$--bridge knots to obtain:

\begin{theorem}[Kobayashi~\cite {Ko4}] \label{2-bridge}
Let $K$  be a non-trivial $2$--bridge knot. Then, for each $g \geq
3$, every genus $g$ Heegaard splitting of the exterior $E(K)$ of $K$
is stabilized.
\end{theorem}

Now using a theorem of Hagiwara \cite{Ha}, stating that all genus three
Heegaard splittings of $E(K)$ obtained from the six genus two ones by
stabilization are isotopic, we obtain:

\begin{theorem}[Kobayashi~\cite{Ko4}] \label{mutualyiso}
Let $K$ be a non-trivial  $2$--bridge knot.  Then, for each $g \geq 3$,
the genus $g$ Heegaard splittings of E$(K)$ are mutually isotopic, that is,
there is exactly one isotopy class of Heegaard splittings of genus $g$.
\end{theorem}

Thus Heegaard splittings of $2$--bridge knots are completely understood
and are obtained by stabilizing and destabilizing one from another.
 
\subsubsection{General tunnel number one knots and links}

For general tunnel number one knots there are some results given mostly
in the terminology of restrictions on the knot type, assuming it has a
single unknotting tunnel. We summarize below some of these results.

For genus one knots, that is, knots with a Seifert surface of genus one,
Scharlemann proves the following:

\begin{theorem}[Scharlemann~\cite{Sc}]
Suppose $K \subset  S^3$ has tunnel number one and genus one. Then either
\begin{enumerate}
\item $K$ is a satellite knot or
\item $K$ is a 2--bridge knot. 
\end{enumerate}
\end{theorem}

Tunnel number one satellite knots were classified by Morimoto--Sakuma
in \cite{MS3} as follows:

A knot $K \subset S^3$ is a $K(\alpha,\beta,p,q,)$ knot if it is
constructed by gluing a non-trivial torus knot $ K_0 = T(p,q)$ to a
$2$--bridge link  $L = K_1 \cup K_2 = K(\alpha,\beta)$,  where $\alpha
\geq 4$: The component $K_2$ of $L$ is a trivial knot so there is a
homeomorphism $h\co E(K_2) \rightarrow N(K_0)$ taking a meridian of the knot
$K_2$ to a fiber $f \subset \partial E(K_0)$ of the Seifert fibration of
$E(T(p,q))$.  The resulting knot  $h(K_1)$ in $S^3$ is denoted as above.

Unknotting tunnels for $K(\alpha,\beta,p,q,)$ are obtained as follows:
Let $\tau_1$ and $\tau_2$ be the top and bottom tunnels for $L$. Set $a_i
= \tau_1\cap \partial N(K_i)$ and $b_i = \tau_2\cap \partial N(K_i),  i=
1,2$. We can assume that $a_2$ is the end point on $\partial E(T(p,q))$
of the unknotting tunnels $\tau_u, \tau_v$ of $E(T(p,q))$ (as in
\fullref{torusknots}). For each $i \in \{1,2\}$  let $\eta_i \subset
\partial N(K_i)$ be an arc joining $a_i$ to $b_i$. Set  $\tau(1,u)
= \tau_1 \cup \tau_u$, $\tau(1,v) = \tau_1 \cup \tau_v$, $\tau(2,u)
= \tau_2 \cup \eta_2 \cup \tau_u$ and $\tau(2,v) = \tau_2 \cup \eta_2
\cup \tau_v$.  Let $D$ denote the Dehn twist in $E(K(\alpha,\beta,p,q,))$
along $\partial E(T(p,q))$ in the direction  of a preferred longitude $l_2
\in \partial E(K_2)$. Now choose $i \in \{1,2\}, w \in  \{u,v\}$ and $n
\in \mathbb{Z}$ and set  $\tau(i,w,n) = D^n(\tau(i,w))$. We can now state:

\begin{theorem}[Morimoto--Sakuma \cite{MS3}]
A tunnel number one non-simple knot in $S^3$ is equivalent to
$K(\alpha,\beta,p,q,)$. Any unknotting tunnel for $K(\alpha,\beta,p,q,)$
is isotopic to  $\tau(i,w,n)$and homeomorphic to  $\tau(i,w)= \tau(i,w,0)$
for some  $i \in \{1,2\}, w \in  \{u,v\}$ and $n \in \mathbb{Z}$.
\end{theorem}

\begin{remark} \label{nonisotopic}
Note that the above theorem determines a whole infinite class of
knots, that is, $K(\alpha,\beta,p,q,)$ which have finitely many minimal
genus Heegaard splittings up to homeomorphism but infinitely many such
splittings up to isotopy.

The first such results are due to Sakuma~\cite{Sa} which are not
published. He produces toroidal $3$--manifolds obtained by gluing Seifert
fibered spaces with base space a disk and two exceptional fibers $D^2
(a/b, c/d)$ to  a $2$--bridge knot complement of type $K(p/q)$ so that
a fiber glues to a meridian.  He further requires that $a/b \neq c/d$
and that  $q^2 \neq 1 \,mod\, p$.  These manifolds are by work of
Kobayashi~\cite{Ko6} of genus two. However one can obtain different
genus two Heegaard surfaces by Dehn twisting along the gluing torus.
Then using the fact that genus two Heegaard surfaces have an involution
and information about the mapping class group he shows that they are
non-isotopic but homeomorphic.
\end{remark}

\begin{remark}
The generalized Waldhausen Conjecture  says that a closed, orientable and
atoroidal $3$--manifold has only finitely many Heegaard splittings of any
fixed genus, up to isotopy.  This conjecture was proved by
Johannson~\cite{Jo,Jo1} for Haken manifolds, which include knot
exteriors. Hence
the above theorem of Sakuma and Morimoto is the best possible.
\end{remark}

We also mention some results for tunnel number one links:

An $n$--string $ n \geq 1$ tangle $T$ is a pair $(B, s)$, where $B$ is a
$3$--ball and $s$ is a finite collection of disjoint  simple closed curves
and $n$ properly embedded arcs. The tangle $T$ is essential if $\partial B
- int(N(s))$ is incompressible and boundary incompressible. A knot or link
$L \subset S^3$ is an {\it $n$--string composite} if the pair $(S^3,L)$
has a decomposition into two essential $n$--string tangles. It will be
called a {\it tangle composite} if it is an $n$--string composite for
some $n$. If it is not a $n$--string composite it will be called  {\it
$n$--string prime}.

\begin{theorem}[Morimoto \cite{Mo}]
A composite link in $S^3$ is tunnel number one if and only if it is a
connected sum of a $2$--bridge knot and the Hopf link.
\end{theorem} 

\begin{theorem}[Gordon--Reid \cite{GR}]
A tunnel number one knot is $n$--string prime for all $n$.
\end{theorem}

\begin{theorem}[Gordon--Reid \cite{GR}]
A tunnel number one link in $S^3$ is tangle composite if and only if it
has an $n$--string Hopf tangle summand for some $n$.
\end{theorem}

\subsection{Knots in $2n$--plat projections}

Historically the first examples of inequivalent Heegaard splittings
for exteriors of knots in $S^3$ were discovered for tunnel number
one knots.  In this sub-section we discuss the method for dealing with
inequivalent Heegaard splittings of higher genus. The method is related
to \fullref{rank} above. For details see Lustig and Moriah \cite{LM2}.

We recall \fullref{NE}: Let $\bar x = \{x_1, \ldots, x_n\}$ and $\bar y
= \{y_1, \ldots, y_n\}$ be two systems of generators for a group $G$.
Denote by $F(X)$ and $F(Y)$ the free groups on bases  $X = \{X_1,
\ldots, X_n\}$ and $Y = \{Y_1, \ldots, Y_n\}$ respectively. Let $\beta_x$
and $\beta_y$ be  the canonical epimorphisms $F(X) \rightarrow G$
given by $X_i \rightarrow x_i$ and $F(Y) \rightarrow G$ given by $Y_i
\rightarrow y_i$. If the generating systems $\bar x$  and $\bar y$ are
{\it Nielsen equivalent} then there is an isomorphism $\alpha\co F(Y)
\rightarrow F(X)$ such that $\beta_y = \beta_x \circ \alpha $. This means
that the words $\{\alpha(Y_1), \dots, \alpha(Y_n)\} \in F(X)$ are a set of
generators and hence are related to $X = \{X_1, \ldots, X_n\}$ by a sequence
of Nielsen automorphisms. These are the collection over all  $i,j \in \{1,
\dots, n\}$ of   the automorphisms given by $X_i \rightarrow X_i^{-1}$
and for  $ i \neq j$, $ X_i \rightarrow X_iX_j$ and their inverses. Given
another system of generators $W = \{W_1, \ldots, W_n\}$ for $F(X)$, written
as words in $X$, we can use the Fox derivative in  $F(X)$ to obtain a
``Jacobian" matrix  $[ \partial W_i /\partial X_j]_{i,j}$. The Jacobian
matrix  of a system of generators obtained by a single application of a
Nielsen automorphism with respect to the original system is an elementary
matrix.  Recall that the Fox derivative satisfies the chain rule. Hence
the matrix
$$[ \partial\alpha(Y_i)/\partial(X_j)]_{i,j}$$
is a product of elementary matrices. The problem, of course, is that we
do not know the automorphism $\alpha$. For any given $y_i \in \bar y$ the
element $\alpha(Y_i)$ is a lift of $y_i$ to $F(X)$ by $\beta_x$.  It is a
theorem that any two lifts to $F(X)$ of $y _i$ differ by some element in
the two sided ideal $I_x \subset \mathbb{Z} F(X)$ generated by $[ \partial
R_i/\partial(X_j)]_{i,j}$ where $R_i$ is one of the relators. Thus
the matrix $[ \beta_x(\partial\alpha(Y_i)/\partial(X_j))]_{i,j}  $
has  entries in $\mathbb{Z} (G) / I_x$,  is  well  defined  there and
is independent of $\alpha$. Furthermore  it  is a product of elementary
matrices.

If there is a non trivial representation $\sigma\co  \mathbb{Z} (G)/
I_{\hat X} \rightarrow \mathbb{M}_m(R)$ then, since $R$ is commutative,
we can compute the determinant of the image matrix over $R$.  If the two
systems $\bar x$ and $\bar y$ are Nielsen equivalent then the determinant
must be a unit of $R$. Showing that it is not so will prove that $\bar
x$ and $\bar y$ are not Nielsen equivalent.  As a  consequence we have
that if  $\bar x$ and $\bar y$ are determined by Heegaard splitting then
these Heegaard splittings are not isotopic.

Using this method  for $m = 1$ and $ R = \mathbb{Z} [t, t^{-1}]$ the
following theorems were obtained:

\begin{theorem}[Lustig--Moriah \cite{LM3}] \label{0.1}
Given a generalized Montesinos knot/link $K \subset S^3$ in a $2n$--plat
projection with associated invariants $\alpha = \gcd \{\alpha_{i,j}\}$
and $\beta = \prod_{i,j}  \{\beta_{i,j}\}$:
\begin{enumerate}
\item If $\beta \neq \pm1$ mod $\alpha$ then the top and bottom tunnel
systems $\tau_t$ and $\tau_b$ are non-isotopic.

\item If $M$ is the closed $3$--manifold obtained by $p_i/q_i$ surgery
on the components $K_i$ of $K$ with $p_i$ even, then $\alpha \neq 1$
implies that $\rank(\pi_1(M)) = \genus (M) = n$.  If $K$ is a knot and
$\beta \neq \pm1$ mod $\alpha$ and $2^{n-1} \neq \pm1$ mod $\alpha$
then the Heegaard splittings  of $M$ induced by $\tau_t$ and $\tau_b$
are non-isotopic.
\end{enumerate}
\end{theorem}

To obtain the full strength of the theory we need to make the following
definitions:

Let $L$ be the $n{\times}2n$--matrix $(a_{i,j})$ with
$$(a_{i,j}) = \begin{cases} 1 & \text{for } j = 2i - 1 \\[-1ex]
  -1 & \text{for } j = 2i \\[-1ex] 0 & \text{otherwise}  \end{cases}$$
Let $M$  be the $2n{\times}n$--matrix $(a_{i,j})$ with
$$(a_{i,j}) = \begin{cases} 1 & \text{for } i = 2j - 1 \\[-1ex]
  1 & \text{for } i = 2j  \\[-1ex]
  0 & \text{otherwise}  \end{cases}$$
Let $N$ be the $n{\times}2n$--matrix $(a_{i,j})$ with
$$(a_{i,j}) = \begin{cases} 1 & \text{for } j = 2i - 1 \\[-1ex]
  0 & \rm{otherwise}  \end{cases}$$
With these conventions we have the following  theorem for a very general
link $K \subset S^3$ with a $2n$--plat projection:

\begin{theorem}[Lustig--Moriah~\cite{LM3}] \label{general}
Let $B$ denote the $2n$--braid underlying the $2n$--plat projection
of $K$. Consider the image $\hat\rho(B)$ of $B$ under the Burau
representation of the $2n$--braid group, with the variable $t$ evaluated
at $-1$. Let $\alpha$ be the greatest common divisor of the entries  of the matrix $L
\circ \hat\rho(B) \circ M$. Let $\beta$ be the determinant of the matrix
$N \circ \hat\rho(B) \circ M$. Then

\begin{enumerate}
\item If $\alpha \neq 1$  then $\rank(\pi_1(E(K))) =  g(E(K))  =  b(K) = n$.

\item If $K_1$ and $K_2$ are two knots, $K = K_1\# K_2$ and  $\alpha
\neq 1$ then $t(K) = t(K_1) + t(K_2)$.

\item If $\alpha \neq 1$ and $\beta \neq \pm1$ mod $\alpha$ then the top
and bottom tunnel systems $\tau_t$ and $\tau_b$ for $K$ are non-isotopic.

\item If $M$ is the closed $3$--manifold obtained by $p_i/q_i$ surgery
on the components $K_i$ of $K$ with $p_i$ even, then $\alpha \neq 1$
implies that $\rank(\pi_1(M)) = \genus (M) = n$.  If $K$ is a knot and
$\beta \neq \pm1$ mod $\alpha$ and $2^{n-1} \neq \pm1$ mod $\alpha$
then the Heegaard splittings  of $M$ induced by $\tau_t$ and $\tau_b$
are non-isotopic.
\end{enumerate}
\end{theorem}

That is, the statements of \fullref{montesinos}  and statements $(1)$
and $(2)$ of \fullref{0.1} hold for $K, \alpha$ and $\beta$.

\begin{remark}
It should be emphasized  that the class of knot/links which are described
by the condition of \fullref{general} is much larger than that of
generalized Montesinos knots/links.
\end{remark}

We end this section with a result obtained by the method discussed
above to distinguish minimal genus Heegaard splittings. It shows that
the Heegaard structure of  minimal genus Heegaard splittings for knot
exterior spaces is indeed a rich one.

\begin{theorem}[Lustig--Moriah \cite{LM4}]
For each $g \geq 3$ there are infinitely many hyperbolic $3$--manifolds
of genus $g$ with at least $2^g - 2$ pairwise non-homeomorphic minimal
genus Heegaard splittings. These are obtained by $m/n$ surgery, $m$
even, on a Montesinos knot $K = {\bf m}(e; (\alpha_1,\beta_1), \dots,
(\alpha_g,\beta_g))$ where the $\beta_i$'s  are mutually distinct
odd primes and $\alpha = \gcd(\alpha_1, \dots, \alpha_g) > 2^{2g -
1}(\beta_1\cdot \dots \cdot \beta_g)^2$.  The exterior $E(K)$ has at
least $2^g - 1$ mutually non-homeomorphic unknotting tunnel systems, that
is,
$2^g - 1$ mutually non-homeomorphic minimal genus Heegaard splittings.
\end{theorem}

\section{Strongly irreducible and weakly reducible}
\label{Strong+weak}

As all knot exteriors are irreducible, a given Heegaard splitting for
a knot space is either stabilized or  if not it is either strongly
irreducible or weakly reducible as in \fullref{reducibility}. So the
question of determining the nature of a given Heegaard splitting is
split into three:

\begin{enumerate}

\item  How do we show that a Heegaard splitting for $E(K)$ is stabilized?

\item  How do we show that a Heegaard splitting of $E(K)$ is  strongly
irreducible?

\item  How do we show that a Heegaard splitting of $E(K)$ is weakly
reducible but not stabilized?
\end{enumerate}

Question (1) was discussed in \fullref{genus}. If the Heegaard splitting
is of minimal genus then it is clearly non-stabilized. If it is not
minimal genus then the known algebraic techniques at our disposal  fail
and in some special cases, that is, exteriors of $2$--bridge knot some things
can be said (see \fullref{2-bridge} and Kobayashi \cite{Ko4}).

If the Heegaard splitting is strongly irreducible then it is clearly
non-stabilized and then the distinction between minimal and non-minimal
genus is not relevant to answering Question (2) above. Below we will
describe some results regarding Question (2).

If it is weakly reducible then it could be stabilized or not depending on
whether the distance of the Heegaard splitting, as in \fullref{distance},
is zero or one. If the Heegaard splitting is of minimal genus then the
fact that it is non-stabilized is obvious and below we will show  examples
of such Heegaard splittings. We are left with  weakly reducible non
minimal genus Heegaard splittings.  It turns out that to decide whether
a Heegaard splitting is distance one if it is not minimal is a very hard
question indeed. In fact the  following remark should be emphasized:

\begin{remark} \label{nemesis}
There is no known example of a non-minimal genus distance one Heegaard
splitting of a $3$--manifold which is either closed or with a single
boundary component\footnote{For further discussion of this issue see
discussion at the end of \fullref{nonmin}.}.
\end{remark}

\subsection{Strongly irreducible Heegaard splittings of knot spaces}
\label{strhs}

We begin with a series of definitions leading to the definition of the
Rectangle Condition due to Casson and Gordon.

\begin{definition}
Let $S$ be a genus $g$ orientable surface and $\mathcal C = \{ c_1, \dots,
c_k\}$ a collection of disjoint essential simple closed curves on $S$.
A {\it wave} with respect to $\mathcal C$ is an arc $\omega$ so that:

\begin{enumerate}
\item  $\partial \omega = \{p_1, p_2\}$ is contained in a single component
$c_{j_0} \in \mathcal C$. 

\item $\omega$ meets $c_{j_0}$ from the same side.

\item  $\omega$ is not homotopic  rel. $\partial \omega$ into $c_{j_0}$.
\end{enumerate}
\end{definition}

Let  $V$ be a  handlebody/compression body and let $ \mathcal C = \{ c_1,
\dots, c_{k}\}$ be a maximal collection of disjoint non-isotopic essential
simple closed curves on $\partial V$  which bound disks $\mathcal D = \{
d_1, \dots, d_{k}\} \subset V$ (that is, $\mathcal D$ is a simplex of maximal
dimension in $\mathcal D_V$ as in \fullref{diskcomplexes}). The boundary
$\partial D$ of any essential disk $D \subset V$ is either parallel to
one of the disks in $\mathcal D$ or has at least two waves with respect
to $\mathcal C$. Such a collection $\mathcal D$ will be called a {\it
pair of pants decomposition for $V$}. A $3$--ball component of $V -
\mathcal D$ will be called a {\it solid pair of pants}.

\begin{definition}[Rectangle Condition]\label{rc}
A Heegaard splitting $(V,W)$ for a  $3$--manifold $M$ will  satisfy the
{\it Rectangle Condition} if there is a {\it blue}  pants decomposition
$\mathcal D$ for $V$ and a a {\it red}  pants decomposition $\mathcal E$
for $W$ so that every pair of  blue curves in every solid pair of pants
of $\partial \mathcal D$ meets every pair of  red curves of every solid
pair of pants in $\partial \mathcal E$.
\end{definition}

\begin{remark}
If $V$ is a compression body with $\partial_-V = T^2$ we assume that a
maximal collection of disks in $V$ will include a separating disk $D_0$
that cuts off a $T^2 \times I$  component from $V$. In this case a
maximal collection $\mathcal D$ will contain $3g - 4$ disks. Note that
no essential disk in $V$ can have a wave in $T^2 \times I$.
\end{remark}

\begin{remark}
Note that a Heegaard splitting which satisfies the Rectangle Condition is
strongly irreducible: Any disk in $V$ ($W$) is either parallel into one
of the curves $d_i$ ($e_j$) or has a wave with respect to $\mathcal D$
($\mathcal E$). A wave travels parallel to the two disks in the pair
of pants which contains it but do not contain the wave's endpoints.
Hence the Rectangle Condition implies that any two essential disks
intersect. In fact they intersect in waves.
\end{remark}

Note also that the Rectangle Condition is strictly stronger than that of
strongly irreducible: The standard Type 2 splittings of $(\text{surface})
{\times}I$, are strongly irreducible but have  pairs of essential disks
that intersect exactly twice.    However  Sedgwick proves in ~\cite{Se2}:

\begin{lemma}
Let $V\cup_S W$ be a Heegaard splitting. Let $D \subset V, E \subset W$
be essential disks. If there are complete collections of disks $\Delta
\subset V,  \Delta' \subset W$ which satisfy the  Rectangle Condition
then $|\partial D \cap \partial E| \geq 4$.
\end{lemma}

\begin{remark}
It should also be pointed out that there are stronger notions then
the Rectangle Condition. This is the {\it Double Rectangle Condition}
of Lustig--Moriah see ~\cite{LM7}.
\end{remark}

\begin{remark}
Given a Heegaard splitting for a $3$--manifold or a knot space with a
pants decomposition, the question of, does it satisfy the rectangle
condition, becomes a combinatorial one. So one can view having the
Rectangle Condition property as a combinatorial  approximation to that
of being strongly irreducible.
\end{remark}

\begin{remark} \label{star}
Casson and Gordon who defined the Rectangle Condition in
unpublished work (see Moriah and Schultens~\cite{MS1}) used it to show
that there are closed $3$--manifolds which have strongly irreducible
Heegaard splittings. A detailed description of these examples can be
found in work of Sedgwick ~\cite{Se}.
\end{remark}

\begin{theorem}[Casson--Gordon]
Let $K = K(p_1, \dots, p_r)$ be a pretzel knot in $S^3$ with  $r  >
5$ strands,  $p_i > 2$ and each $p_i$ odd. Let $M$ be a manifold which
is  obtain by  $\frac{1}{k}$--surgery  for $k \geq 6$. Then the Heegaard
splitting of $M$ obtained from a regular neighborhood of a Seifert
surface for $K$ and its complement is strongly irreducible.
\end{theorem}

In \cite{LM6} M Lustig and the author also used the Rectangle Condition
to show that in fact the  phenomenon of manifolds with strongly
irreducible Heegaard splittings is quite common and that there are many
knots in $S^3$ which produce after surgery such $3$--manifolds.  Kobayashi
in ~\cite{Ko1} also exhibits manifolds with strongly irreducible Heegaard
splittings. These examples will be further discussed in \fullref{nonmin},
\fullref{infHS}. Recently Minsky, Moriah and Schleimer have showed
in ~\cite{MMS} that for any positive integer  $t \geq 1$ and $n > 0$
there are tunnel number $t$ knots   $K \subset S^3$ so that $E(K)$ has a
Heegaard splitting of distance greater than $n$. In particular  they are
all strongly irreducible and  in fact it seems knots with high distance
Heegaard splittings are the ``generic" case.

Many of the definitions and results which come under the above ``Strongly
irreducible Heegaard splittings of knot spaces" heading of  this section
are also discussed  in depth  in Sections \ref{infHS} and  \ref{super}.
We refer the reader to those sections.

\subsection{Weakly reducible non-stabilized Heegaard splittings} \label{wrhs}

The previous section tells us that knots in $S^3$ with strongly
irreducible Heegaard splittings are abundant. So where does one look
for knots with weakly reducible Heegaard splittings?  The first  such
result is due to M Lustig and the author ~\cite{LM5}. In that
paper the notion of a {\it wide knot} is defined as follows:

\begin{definition}
\begin{itemize}
\item[(a)] A $2n$--braid will be called {\it wide} if in its standard
projection where every crossing is replaced by a node, there is no
monotonically descending path connecting the top of the second strand
to the bottom of the ($2n-1$)-strand or vice versa.

\item[(b)]  A $2n$--plat projection of a knot or link will be called
{\it wide} if the underlying braid is wide.

\item[(c)] A knot or link $K \subset S^3$ will be called {\it wide}
if it has a wide $2n$--plat projection so that the corresponding Heegaard
splitting is irreducible.
\end{itemize}
\end{definition}

Then the following proposition is proved:

\begin{proposition}[Lustig--Moriah~\cite{LM5}]
For every knot or link $K$ in $S^3$ in a wide $2n$--plat projection both
canonical Heegaard splittings are weakly reducible. In particular every
wide knot or link has a weakly reducible  and irreducible (distance one)
Heegaard splitting.
\end{proposition}

\textbf{Idea of Proof}\qua Consider the canonical Heegaard splitting
of $E(K)$ induced by inserting the bottom tunnels to the $2n$--plat
projection. Embed a $2$--sphere $S$ in $(S^3,K)$ which cuts off the ``top"
bridges. When a regular neighborhood $N(K)$ is removed from $S^3$ to
obtain $E(K)$ the ``inside" of $S$ (that is, the $3$--ball (before removing
$N(K)$) containing the top bridges) is isotopic to $W$ of the Heegaard
splitting. It is clear that the leftmost bridge  defines a dual essential
disk $D$ in $W$. Since the knot/link is wide then in the process of
gluing the handlebody to the compression body the boundary of the disk
$D$ will be mapped according to the braid. Hence it will not meet the
cocore disk of the rightmost bottom tunnel which is an essential disk
in the compression body component $V$ of the Heegaard splitting.

The fact that the Heegaard splitting is irreducible is proved for many
knots/links in a $2n$--plat projection using the algebraic techniques
of \fullref{montesinos}.

Another place to look for knots which have Heegaard splittings of
distance one are knots which are connected sums.  This approach was
taken by the author in ~\cite{Mr2}. We first need some definitions.

\begin{definition}\label{primitive}
We say that a curve on a handlebody is {\it primitive} if there is
an essential disk in the handlebody intersecting the curve in a single
point. An annulus $A$ on the boundary of a handlebody is  primitive if its
core curve is a primitive curve. A Heegaard splitting $(V, W)$ for $E(K)$,
where $\partial E(K) \subset V$ will be  called {\it $\gamma$--primitive}
if there is a spanning annulus $A\subset V$ such that  $\partial A =
\gamma \cup \alpha$ where $\gamma \subset \partial E(K)$ and $\alpha$
is a primitive curve on the Heegaard surface $\partial W$.  If $\gamma$
is a meridian for $K$ then we say that $(V,W)$ is {\it $\mu$--primitive}.

We say that a knot $K \subset S^3$ has a primitive meridian, or is {\it
$\mu$--primitive} if $E(K)$ has a minimal genus Heegaard splitting with
a primitive meridian.
\end{definition}

\begin{definition}\label{induced}
Two Heegaard splittings $(V_i, W_i)$ for  $E(K_i)$ respectively, induce
a decomposition of $E(K) = E(K_1 \# K_2)$ into $(V, W)$. We can think
of $V_i$ as a union of  $\partial E(K_i) \times I$ and $1$--handles. Let
$t_i \subset K_i$ be an arc.  Let  $B_i$, a small $3$--ball in $S^3$
such that $t_i$ is an unknotted properly embedded arc in $B_i$, the
closure of $B_i \cap S^3 - E(K_i)$ is a regular neighborhood $N(t_i)$
of $t_i$ in $B_i$ and the annulus $A'_i = \partial B_i - N(\partial t_i)$
is the union of two vertical annuli  $A_1^{*i}, A_2^{*i}$ and a meridional
annulus $A_i\subset \partial E(K_i) \times \{1\}  \subset  \partial_+V_i =
\partial W_i$. We regard $(B_i, t_i)$ as the pair which is removed from
$(S^3, K_i)$ when forming the connected sum of pairs $(S^3, K_1 \# K_2) =
(S^3, K_1)\#(S^3, K_2)$.  We obtain $V$ as follows:

Cut the compression bodies $V_1$ and $V_2$ open along the annuli
$A_1^{*1}, A_2^{*1}$ and $A_1^{*2}, A_2^{*2}$ respectively.  Consider
the components which are not  $B_i - N(t_i)$.  We obtain handlebodies
$V'_1$ and $V'_2$ with  copies of  $A_1^{*1}, A_2^{*1} \subset V'_1$
and copies of $A_1^{*2},A_2^{*2} \subset V'_2$. Now glue $V'_1$ to $V'_2$
by identifying $A_1^{*1}$ with $A_1^{*2}$ and $A_2^{*1}$ with $A_2^{*2}$
to get  a compression body $V$. We obtain $W$ by gluing $W_1$ and $W_2$
along the meridional annulus $A_i$.  Hence $V$ is always a compression
body but $W$ is a handlebody if and only if the  meridional annulus $A_i$
is a primitive annulus in $\partial W_i$ for one of $i = 1 $ or $i = 2$.
In this case we  will say that $(V, W)$ is the {\it induced Heegaard
splitting of  $E(K)$} induced  by  $(V_1, W_1)$ and $(V_2, W_2)$.
\end{definition}

\begin{theorem}[Moriah ~\cite{Mr2}] \label{wredProp}
Let $K_1, K_2$  and  $ K = K_1 \# K_2$ be knots  in $S^3$ and $(V_i ,
W_i ), i =1,2$ be  Heegaard splittings  for $E(K_i)$. If $(V_1, W_1)$
and $(V_2 , W_2)$ induce a Heegaard splitting $(V, W)$  of  $E(K)$
then $(V, W)$ is a weakly reducible Heegaard splitting.
\end{theorem}

Using the above definitions we have the following two theorems:

\begin{theorem}[Moriah ~\cite{Mr2}] \label{Supaddthm}
Given knots $K_1, K_2$  and  $ K = K_1 \# K_2$ in $S^3$  for which the
tunnel number satisfies $t(K) = t(K_1) + t(K_2) +1$  then there is an
induced  minimal  genus Heegaard splitting of  $E(K)$ which is  weakly
reducible.
\end{theorem}

These theorems raise the following two questions:

\begin{question} \label{strmin}
Are there non-prime knots in $S^3$ so their exteriors have  strongly
irreducible minimal genus Heegaard splittings?
\end{question}

\begin{question} \label{both}
Are there knots in $S^3$ so their exteriors  have
both weakly reducible and strongly irreducible minimal genus Heegaard
splittings?
\end{question}

\begin{remark} \label{down}
There is a positive ansewer to \fullref{strmin}. The first example
was given by the author in ~\cite{Mr2}. There  $K_1= K^n(-2,3,-3,2)$
is the twisted pretzel where between the two $3$--strands of twists one
introduces  an odd number $n \in \mathbb{Z},  n \notin \{-1, 0, 1\}$ of
``horizontal" crossings, and $K_2$ is any $2$--bridge knot. Then  $t(K_1)
= 2$  and $t(K_2) = 1$. It has been used by Morimoto ~\cite{Mo1}
to show that the tunnel number can go down after connected sum, that is, it
is sub-additive. He exhibits a Heegaard splitting of genus three so that
$t(K_1 \# K_2) = 2 < 2 + 1 = t(K_1) + t(K_2)$ (see \fullref{additivity},
\fullref{sub}). We have:
\end{remark}

\begin{proposition}[Moriah ~\cite{Mr2}]
Let $K_1$ and $K_2$ be as above. The genus three Heegaard splitting
of $E(K_1 \# K_2)$  described in ~\cite{Mo1} is strongly irreducible.
\end{proposition}

This result was generalized by Kobayashi and Rieck in ~\cite{KR}. They prove:

\begin{theorem} [Kobayashi and Rieck {{\cite[Corollary 5.4]{KR}}}]
Let $K_1$ and $K_2$ be knots in closed orientable 3--manifolds. Let $X_i =
E(K_i) (i = 1, 2)$ and $X = E(K_1\#K_2)$. Suppose that:
\begin{enumerate}
\item  X is irreducible, 
\item $X_i,  (i = 1, 2)$ does not contain a meridional essential annulus, 
\item $X_i , (i = 1, 2)$ does not contain an essential torus, and
\item $g(X_1) + g(X_2) \geq 5$ and $g(X) = 3$.
\end{enumerate}
Then any minimal genus Heegaard splitting of $X$ is strongly irreducible. 
\end{theorem}

\fullref{both} is still open for manifolds which are exterior of knots in
$S^3$. However for general $3$--manifolds it has been answered positively
by Kobayashi and  Rieck. They have proved:

\begin{theorem}[Kobayashi and Rieck  ~\cite{KR}]
There are infinitely many $3$--manifolds which have both strongly
irreducible and weakly  reducible Heegaard splittings of minimal genus.
\end{theorem}

\textbf{Idea of proof}\qua The manifolds $M$ are constructed as follows:
Let $X = S^3 - N(K_1)$, $Y = S^3 - N(K_2)$ and  $Z = S^3 - N(K_3)$,
where $K_1 = T(2,3)$ is the trefoil knot, $K_2 = L(\alpha, \beta)$
with $ \alpha $ even, is any $2$--bridge link which is not the Hopf link
and $K_3 = K(2,5) $ is the figure $8$ knot.  Let $\mu_1$ and $\mu_2$
be meridians of  the $2$--bridge link $L(\alpha, \beta)$ components,
$\lambda$ be the longitude of $K_3$ and $\gamma$  be the boundary of
the cabling annulus in $X$.

Attach $X$ and $Z$ to $Y$ by gluing their tori boundaries so that $\gamma$
is mapped to $\mu_1$  and $\lambda$ to $\mu_2$.  We obtain a closed
$3$--manifold $M$.  By analyzing the incompressible tori contained in $M$
they conclude that $g(M) = 3$ and that $M$ has a strongly irreducible
minimal genus  Heegaard splitting.  On the other hand $M$ has a Heegaard
splitting which is the amalgamation of two genus two Heegaard splittings
along a torus. Hence it is of genus three and is clearly weakly reducible.

\section{Non-minmal genus Heegaard splittings}
\label{nonmin}

In general non-minimal genus Heegaard splittings of exteriors of knots in
$S^3$ are not well understood at all. The only comprehensive results are
those of Kobayashi and Moriah--Schultens mentioned in \fullref{t(K)=1},
Subsections \ref{$2$--bridge knots} and  \ref{torusknots}. That is
\fullref{2-bridge} and \fullref{mutualyiso}, where  all Heegaard
splittings of $E(K)$ for $K$ a $2$--bridge knot are classified. Also
complements of torus knots are Seifert fibered spaces and as such
Moriah--Schultens ~\cite[Theorem~0.1]{MS1} states that any of their non-minimal Heegaard
splittings is a stabilization of minimal genus two vertical Heegaard
splittings.

For manifolds with non-empty boundary there is  one other technique
for obtaining non-minimal genus Heegaard splittings other than
stabilizing. It is called ``boundary stabilization", defined by the author
in ~\cite{Mr1}.   Boundary stabilization uses an operation on Heegaard
splittings called ``amalgamation".  Amalgamated Heegaard splittings were
defined by Schultens ~\cite{Sch}. We give the full definition for the
benefit of the reader:

\begin{definition}\label{amalgamation}
Given two manifolds $M_1$ and $M_2$ with respective Heegaard splittings
$(U^1_1,U^2_1)$ and $(U^1_2,U^2_2)$, assume further that there are
homeomorphic boundary components $F_1 \subset \partial_{-} U^1_1$
and  $F_2 \subset \partial_{-} U^1_2$. Denote the homeomorphism $F_1
\rightarrow F_2$ by $g$.  Let $M$ be a manifold obtained by gluing $F_1$
and $F_2$ along the homeomorphism $g$.  We can obtain a Heegaard splitting
$(V_1,V_2)$ for $M$ by a process called {\it amalgamation} as follows:

Given a compression body $U$ we can assume that it has the structure
of $\partial_{-} U \times I  \cup(\cup \{1$--handles$\})$. Let $h_i$
be the homeomorphism $N(\partial_{-} U^1_i) \rightarrow (\partial_{-}
U^1_i) \times I$ and  $p_i\co (\partial_{-} U^1_i) \times I \rightarrow
\partial_{-} U^1_i$ the projection into the first factor.  Define an
equivalence relation  $\equiv$  on $M_1 , M_2$  as follows:

\begin{enumerate}
\item  If $x_i, y_i \in N(F_i)$  are points such that $p_i h_i(x_i) =
p_i h_i(y_i)$  then  $x_i  \equiv y_i$.
\item  If $ x \in F_1 , y \in F_2$  and  $g(x) = y$, where  $g\co F_1
\rightarrow F_2$  is the homeomorphism between the surfaces, then  $x
\equiv y$.
\end{enumerate}

Furthermore we can arrange that the attaching disks on $F_1 \times
\{1\}, (F_2 \times \{1\})$ for the $1$--handles in  $ U^1_1$, $U^1_2$
respectively, have disjoint images in  $F_1$ ($F_2$)  and hence they do
not get identified to each other.  Now set $M = (M_1 \cup M_2) / \equiv $,
$V_1 = (U_1^1 \cup U_2^2) / \equiv$, $V_2 = (U_1^2 \cup U_2^1) / \equiv$.
Note that $V_1 = U_2^2 \cup N(F_1) \cup (\cup \{1$--handles$\}) \cup
(\partial_{-}  U_1^1 - F_1) \times I$ and $V_2 = U_1^2 \cup N(F_2) \cup(
\cup \{1$--handles$\}) \cup (\partial_{-} U_2^1 -F_1) \times I$ so that
$V_1$ and $V_2$  are compression bodies defining a  Heegaard  splitting
$(V_1,V_2)$ for  $M$ .  The  Heegaard splitting  $(V_1,V_2)$  of  $M$  is
called the {\it amalgamation} of the Heegaard splittings $(U_1^1,U_1^2)$
of  $M_1$  and  $(U_2^1,U_2^2)$ of  $M_2$  along  $F_1$ and  $F_2$ .
\end{definition}

We give here a modified version for knot spaces:

\begin{definition}
Given a genus $g$ Heegaard splitting for a knot exterior $E(K)$ in
$S^3$  we can obtain a genus $g + 1$ Heegaard splitting by taking the
standard genus two Heegaard splitting for  $S^1\times S^1 \times I$
(see Scharlemann and Thompson~\cite{ST}) and amalgamating the two splittings via the identity
map $T^2 \times \{1\} \rightarrow \partial E(K)$. The amalgamated Heegaard
splitting will be called a {\it boundary stabilization}.
\end{definition}

It is  clear that if the knot is  $\gamma$--primitive (see
\fullref{primitive} in \fullref{Strong+weak}) then a boundary
stabilized  Heegaard splitting is in fact stabilized: Since the knot is
$\gamma$--primitive then after gluing $T^2 \times I$ to $\partial E(K)$
there is an annulus $A$ in $T^2 \times I  \cup E(K)$ which meets a disk
$D$ in the handlebody component of the non-stabilized Heegaard splitting
in a single point. It is just the extension of the annulus in $E(K)$ which
has the same property. The amalgamation process removes a neighborhood
of an essential arc from $A$ thus creating an essential disk $E$ in the
handlebody component of the boundary stabilized Heegaard splitting.
The disk $E$ still meets the disk $D$ (now in the compression body
component of the boundary stabilized Heegaard splitting) in a single
point  (see ~\cite{Mr1}). So the boundary stabilized Heegaard splitting
is stabilized. However if the knot is not $\mu$ or  $\gamma$--primitive
then the following conjecture is proposed:

\begin{conjecture}
Given a knot  $K \subset S^3$ which is not  $\gamma$--primitive then a
boundary stabilization of a minimal genus Heegaard splitting of $E(K)$
is non-stabilized.
\end{conjecture}

\subsection{Manifolds with arbitrarily high genus irreducible Heegaard splittings}
\label{infHS}

Though this sub-section discusses closed manifolds, knot exteriors played
a crucial role in showing that there are $3$--manifolds with irreducible
Heegaard splittings of arbitrarily high genus.

The actual result of Casson and Gordon discussed in \fullref{strhs},
\fullref{star},  is much stronger than the one stated. They used the Rectangle
Condition to  show that there are closed  $3$--manifolds which have
strongly irreducible Heegaard splittings of, in fact, arbitrarily high
genus.  This result is unpublished and a  proof of it, similar to the
original proof, is given by Kobayashi \cite{Ko}.  A different  proof due to  Casson
is given in the appendix of Moriah--Schultens \cite{MS1}. These manifolds  are obtained by
$\frac 1 k$--surgery  for $k \in \mathbb{Z}, k \geq 6$ on pretzel knots
in $S^3$ of the form $K(p_1, \dots, p_r)$.  These are pretzel knots with
$r  > 5$ strands and each strand has  $p_i > 2$ and odd crossings.

Results of  M Lustig and the author ~\cite{LM6} showed that in fact this
phenomenon is quite common and that there are many knots in $S^3$ which
produce after surgery $3$--manifolds with strongly irreducible Heegaard
splittings of arbitrarily high genus.  These knots are given as knots
which embed on the boundary of regular neighborhoods of planar graphs
$N(\Gamma)$ called {\it trellises}. Trellises generalize the notion of
a $2n$--plat and allow us to  present a knot or link carried by them by
a family of integer  parameters, assembled in a  {\it twist matrix}. We
can compute a {\it trellis linking number a(K)}. It is clear that both
$N(\Gamma)$ and  $S^3 - N(\Gamma)$ are handlebodies. Then one obtains
a closed manifold by doing a  $(1 + k a(K)) / k$--surgery on $E(K)$,
for any  $k \in \mathbb{Z}$. This surgery is so designed to actually
correspond to a Dehn twist on $\partial N(\Gamma)$. A Heegaard splitting
that is created in this way is called {\it horizontal}. Once the two
handlebodies are glued together with this Dehn twist we get the  required
closed  manifold.  The knots satisfy some complexity to ensure that the
Rectangle Condition for these Heegaard splittings will be satisfied so
that we get a strongly irreducible Heegaard splitting. We then have:

\begin{theorem} \label{one}
 Let $K$ be a knot given as a $2n$--plat in  $S^3$, and assume that all
 twist coefficients
of the underlying braid satisfy $|a_{i,j}| \geq 3$. Then for all  $k
\in \mathbb{Z}$, with $|k| \geq 6$, the manifold $K((1 + k a(K)) / k)$
has a strongly  irreducible horizontal Heegaard splitting.
\end{theorem}

There are planar graphs $T$ which are a generalization of trellisses. If
we consider generalized trellisses $T$ with a particular  combinatorial
feature, called an {\it interior pair of edges}, we can perform flypes on
these  more  general knots in a way similar to that done by Casson--Gordon
in ~\cite{MS1} for  pretzel knots.  Define the ``genus" of $T$ by $g(T)
= -2(\chi (T))$.

\begin{theorem} \label{two}
Let $T$ be a generalized trellis and let $K = K(A) \subset S^3$
be a  knot carried by $T$ with twist matrix $A$. Assume that all
coefficients $a_{i,j}$ of $A$ satisfy  $|a_{i,j}| \geq 3$ and that
there is an interior pair of edges $(e_{i,j}, e_{i,h})$ of $T$ with
twist coefficients  $|a_{i,j }|, |a_{i,h}| \geq 4$.   Then for all $k,
n \in \mathbb{Z}$, with   $|k| \geq 6$, the manifolds $K((1 + k a(K))
/ k)$ have irreducible Heegaard splittings $\Sigma(n)$ of arbitrarily
large genus $g(T)+2n$, all of which are horizontal.
\end{theorem}

\begin{remark}
In particular, for $K$ as in \fullref{two}, all of the splittings
$\Sigma(n)$, stabilized once, are stabilizations of a common low genus
Heegaard splitting.
\end{remark}

There are other examples of manifolds with strongly irreducible Heegaard
splittings of arbitrarily high genus.  These manifolds were introduced
by Kobayashi in  ~\cite{Ko1}. Kobayashi considers two component links
$L^n = l_1^n \cup l_2^n \subset S^3$ which are pretzels of the form
$$P(9,-9,7,5,-5-7,5,5,-5,-5, \dots,5,5,-5,-5),$$
where there are $n-1$ tangles of the form $(5,5,-5,-5)$. Two Seifert
surfaces $S, S'$ in $E(L)$ will be called {\it weakly equivalent} if
there is a homeomorphism $h\co E(L) \rightarrow E(L)$ so that  $h(S) = S'$.
For these links we have:

\begin{proposition}[Kobayashi~\cite{Ko1}]\label{poly}
For each integer $ g \geq n $,  the link $ L^n $ has
$\binom{g-1}{n-1}$ free, incompressible Seifert surfaces $S$
of genus $2n + g$ which are mutually non-weakly equivalent.
\end{proposition}

\fullref{poly} is proved by showing that different surfaces are carried
by an  incompressible branched surface  with a different set of weights
$\{m_1, \dots m_n\}$ so that $m_1 + m_2 + \dots m_n = g$. Two surfaces are
weakly equivalent if and only if $m _i = m'_i$ for each $i = 1, \ldots, n$.

Now let $L = l_1 \cup l_2$ be a $2$--bridge link and $f\co \partial E(L)
\rightarrow \partial E(L^n)$ be an orientation reversing homeomorphism
mapping a meridian of $l_i$ to a longitude of $l_i^n$.  Define $M = E(L)
\cup_f E(L^n)$. This gluing induces a strongly irreducible Heegaard
splitting of $M$ of genus $4n + 2 +2g$. If $g > n$ these Heegaard
splitting are not of minimal genus.

For the Heegaard splittings of the manifold $M$ we have:

\begin{proposition}[Kobayashi~\cite{Ko1}]\label{weakequiv}
If $L$  is not a $(2,2n)$ torus link and two Heegaard splittings of $M$
induced by the Seifert surfaces $S$ and $S'$ are homeomorphic then the
homeomorphism is a weak equivalence between $S$ and $S'$.
\end{proposition}

Note that in all the examples so far~ (Kobayashi \cite{Ko,Ko1},
Lustig--Moriah \cite{LM6} and Moriah--Schultens \cite{MS1}), the knots/links $K \subset S^3$ in question all have
the property that the higher, non-minimal genus splittings are obtained
by iteratively twisting the knot in a Conway sphere through a $\pi$
angle via an ambient isotopy of $S^3$. Though the knot type stays the
same the new ``position" of the knot induces a new Heegaard  splitting
which has genus bigger by two than the previous one.

One should also note that {\it all}  the Heegaard splittings discussed
above are strongly irreducible and two Heegaard splittings become isotopic
after stabilizing once the higher genus and the corresponding number of
stabilizations for the lower genus  splitting.

A set of somewhat different examples were introduced by Moriah, Schleimer
and Sedgwick in ~\cite{MSS}. These examples are not related to knot
spaces but they are brought here to emphasize
the fact that the phenomenon of manifolds with strongly irreducible
Heegaard splittings of arbitrarily high genus is very common.

Note that the manifolds of Casson--Gordon have Heegaard genus four and
larger.  These new examples have genus as low as three.  Also, these
examples, unlike those of~\cite{Ko1} and~\cite{LM6}, do not involve
twisting around a two-sphere in $S^3$ or require the existence of an
incompressible spanning surface.

Here is a sketch of the construction, which has obvious generalizations:
Let  $V$ be a handlebody of genus three or more.  Choose $\gamma$
to be a ``disk busting'' curve in $S = \partial V$,  that is, a curve
which intersects every essential disk in $V$. Note that $ K' = S -
N(\gamma)$ is an incompressible surface in $V$. Let $W$ be another
copy of $V$ and now double $V$ across $S$ after modifying the gluing
of $V$ to $W$ by Dehn twisting along $\gamma$ at least six times.
This gives a closed orientable manifold  $M$.  As $K' \cap \gamma =
\emptyset$ the surface  $K'$ doubles to give a surface $K$ in $M$. which
intersects $S$ transversally.  Adding copies of $K$ to $S$ via the
Haken sum operation (see ~\cite{MSS}) will give the desired sequence
of Heegaard splittings. This requires  proving the following theorem
(see ~\cite{MSS}):

\begin{theorem}
\label{Thm:NewExample}
Given $V$ and $\gamma$ as above, the surface $S + nK$ is a strongly
irreducible Heegaard splitting of $M$, for any even $n > 0$.
\end{theorem}

It is clear from the construction that the genus of the Heegaard surface
$S + nK$ goes up as $n$ increases.

\begin{remark}
There are examples of complements of links in $S^3$ which have weakly
reducible non-stabilized Heegaard splittings. The idea behind all of
these examples is that each of these spaces has Heegaard splittings which
have different partitions of the boundary components in the compression
bodies. The first such example is due to Sedgwick ~\cite{Se1} and is the
three component unlink in $S^3$. It has a minimal genus two Heegaard
splitting where one compression body contains two of the boundary
components and the other compression body contains one. It also has a
clearly weakly reducible non-stabilized genus three Heegaard splitting
with all three boundary components in  a single compression body. In
further work by Moriah and Sedgwick  ~\cite{MS2} two infinite sets of
different examples of two component links $L \subset S^3$ with genus
three weakly reducible non-stabilized Heegaard splittings of $E(L)$
which are not minimal genus, were found.
\end{remark}

So after all of the results and discussion above we have still made no
progress regarding the nemesis of Heegaard splittings:  Weakly reducible
but irreducible Heegaard splittings?

\section{Additivity of Tunnel number}
\label{additivity}
Given two knots  $K_1, K_2 \subset S^3$ it is a natural question to ask:
How does the genus of $E(K_1)$ and $E(K_2)$ compare with the genus of
the connected sum $E(K) = E(K_1\#K_2)$? This problem is usually referred
to as the {\it additivity of tunnel number}.

It is a well known fact  that gluing two compression bodies $V_1$ and
$V_2$ along an annulus $A$ will result in a compression body $V= V_1
\cup_A V_2$, if and only if $A$ is primitive, as in \fullref{primitive},
in either $V_1$ or $V_2$. For a proof  of this result and an extension
see ~\cite{Mr2}.

Recall that given two knots  $K_1, K_2 \subset S^3$ the exterior $E(K)$
of the connected sum $K = K_1 \# K_2$ is obtained by gluing $E(K_1)$
to $E(K_2)$ along a meridional annulus $A$.

As a consequence, it is a  fact that for all knots $K = K_1\#K_2 \subset
S^3$ the tunnel numbers satisfy $t(K) \leq t(K_1) + t(K_2) +1 $:  Let
$\{\tau_1, \ldots, \tau_{t(K_1)}\}$ be a minimal tunnel system for $K_1$
and $\{\sigma_1, \ldots, \sigma_{t(K_2)}\}$ be a minimal tunnel system
for $K_2$. Let $A$ denote the decomposing annulus for $K = K_1\#K_2$
and  let $a$ be an essential arc on $A$. The arc $a$ is an unknotting
tunnel for say $K_1$ (It can be thought of as an unknotting tunnel for
$K_2$ as well.) As if it is slightly pushed into the handlebody $W_1 =
S^3 - N(K_1 \cup\{\cup \tau_i\})$ in a neighborhood of $A$   it  clearly
has a dual disk in $W_1$. Furthermore $A$ becomes a primitive annulus in
$W_1$ because its core curve intersects the same disk in a single point.
Hence  if $W_2$ is the handlebody complement  of $S^3 - N(K_2 \cup \{\cup
\sigma_i\})$ then  $W_1 \cup_A W_2$  is a handlebody. So $\{a, \tau_1,
\ldots, \tau_{t(K_1)},\sigma_1, \ldots, \sigma_{t(K_2)}\}$ is a tunnel system
of cardinality $t(K_1) + t(K_2) +1 $ for $K$.

In other words the Heegaard splittings of  $E(K_1)$ and $E(K_2)$
determined by
$$\{ a, \tau_1, \ldots, \tau_{t(K_1)}\} \quad\text{and}\quad
  \{\sigma_1, \ldots, \sigma_{t(K_2)}\}$$
induce a Heegaard splitting, as in \fullref{induced}, on $E(K)$.

Hence there are three possibilities for $K$:

\begin{description}
\item[Super additive] $t(K) = t(K_1) + t(K_2) +1 $, that is,
$$g(E(K)) = g(E(K_1)) + g(E(K_2)).$$

\item[Additive] $t(K) = t(K_1) + t(K_2)$, that is,
$$g(E(K)) = g(E(K_1)) + g(E(K_2))  - 1.$$

\item[Sub-additive] $t(K) <  t(K_1) + t(K_2)$, that is,
$$g(E(K)) \leq g(E(K_1)) + g(E(K_2)) - 2.$$
\end{description}

Before we discuss  the three cases we need some definitions and context.

\begin{definition} \label{unknotted}
Suppose $\mathcal T \subset V$ is a disjoint collection of properly
embedded arcs in a compression body $V$ where $\partial \mathcal T \subset
\partial_+ V$.  We say $\mathcal T$ is {\em unknotted } if $\mathcal T$
can be properly isotoped, rel boundary, into $\partial V$.
\end{definition}

The following well-known generalization of bridge position is due to
Doll~\cite{Do}:

\begin{definition}\label{g,b}
Suppose that $M = V \cup_S W$ and $K$ is a knot in $M$.  The knot $K$
is in {\em bridge position with respect to $S$} if $K$ is transverse
to $S$ and either
\begin{enumerate}
\item $K \cap S \neq \emptyset$, and both $K \cap V$ and $K \cap W$
are unknotted, or
\item  $K \cap S = \emptyset$,  $K \subset V$ (without loss of generality)
and $V - N(K)$ is a  compression  body.
\end{enumerate}
If $g{=}g(S)$ and $b{=}|K{\cap}S|/2$ then we say that $K$ admits a {\em
$(g,b)$--decomposition}.  Saying that $K$ has no $(g,b)$--decompositions
means that $K$ admits no $(g',b')$--de\-com\-pos\-i\-tions with $g'\leq g, b'
\leq b$.
\end{definition}

There is a  connection between the notions of additivity of tunnel
number and $(g,b)$--decompositions.  It is exhibited, for example, in
the following known fact:

\begin{lemma}
\label{Lem:Equivalent} 
Let $K \subset M$ be a knot with $t(K) = t$. Then $K$ has a
$(t,1)$--decomposition if and only is $K$ is $\mu$--primitive.
\end{lemma}

We give a proof to illustrate the connection.

\begin{proof} 
Let $(V,W)$ be the Heegaard splitting of genus $t$ of $M$ which realizes
the $(t,1)$--decomposition. Hence $K = t_1 \cup t_2$ where $t_1 \subset V,
t_2 \subset W$ are unknotted arcs. Consider a regular neighborhood of $t_2
\subset W$. We can think of it as a  $1$--handle containing $t_2$. Remove
it from $W$ and add it to $V$ to obtain a handlebody $V'$ of genus $t +
1$. Since $t_2$ was unknotted the drilled out manifold $W'$ is also a
handlebody. The cocore disk of the $1$--handle meets an essential disk of
$W'$ in a single point. This disk is the disk which is determined by the
given isotopy of $t_2$ into $\partial_+W$. When $K \subset V'$ is removed
we get a compression body $V''$. Thus $(V'', W')$ is a Heegaard splitting
of minimal genus $g = t + 1$ and by the construction is $\mu$--primitive.

If $(V,W)$ is a minimal genus Heegaard splitting for $M - N(K)$ with $V$
being the compression body then add a solid torus neighborhood of $K$
to $V$ to obtain a handlebody $V'$. It has an essential disk $D$ which
meets $K$ in a single point. Since $(V,W)$ is $\mu$--primitive there
is an essential disk $E$ of $W$ which $D$ meets in a single point.
If we remove a regular neighborhood $N(D)$ of $D$ from $V'$ and add it
to $W$ along an annulus as a $2$--handle we obtain handlebodies $V''$
and $W''$ and hence a Heegaard splitting $(V'',W'')$ of genus $g -1 =
t(K)$ for $M$. Now $K \cap N(D)$ is an arc $t_2 \subset W''$ and $K
\cap V''$ is an arc $t_1$ and they are clearly unknotted. So we have a
$(t,1)$--decomposition for $K \subset M$.
\end{proof}

We will discuss the topic of additivity of tunnel number according to the
above partition into super additive, additive and sub-additive. However
this partition is somewhat artificial and hard to enforce as various
results in one section are relevant and could be seen as belonging
to another.

\subsection{(1) Super Additive} \label{super}

It clearly follows from the above discussion and Definitions
\ref{primitive} and \ref{induced} that if either of the knots $K_i,
i = 1,2$ is $\mu$--primitive then case (1), as in the trichotomy above,
cannot happen.

\begin{remark}
For a long time it was an open question whether there are any knots that
satisfy super additivity. This was finally resolved by Moriah--Rubinstein
~\cite{MR} and Morimoto, Sakuma and Yokota ~\cite{MSY}.
\end{remark}

In ~\cite{MR} the following theorem was proved:

\begin{theorem}[Moriah and Rubinstein ~\cite{MR}]\label{n+m+1}
For any pair of odd positive integers $r, s \geq 3$ there exist infinitely
many pairs of knots $K_1, K_2 \subset S^3$ so that $t(K_1) = (r -1) /
2$, $t(K_2) = (s - 1) / 2$ and  $t(K_1 \#K_2) = t(K_1) + t(K_2) + 1  =
(r+s) / 2$.
\end{theorem}

The above \fullref{n+m+1} is an existence theorem and is not constructive
in the following sense: The knots $K_1$ and $K_2$ are obtained as the
branched sets in $S^3$ of a $2$--fold cover of $S^3$ by ``sufficiently
large " $m / n$--surgery on pretzel knots $K_1 = K(p_1, \ldots, p_r)$ and $K_2 =
K(p'_1, \ldots, p'_s)$ where $r, s$ are odd, $p_i = p_{r - i +1}$, and $p'_i =
p'_{s - i +1}$. As we cannot determine the precise number $m/n$ we can
only get an existence result.

In ~\cite{MSY}  Morimoto, Sakuma and Yokota actually give the first
concrete examples of knots which are  super additive. These knots are
all tunnel number one knots. They are of the following form:

A knot in $S^3$ will be called a {\it twisted torus knot} of the form
$T(p,q,r,n)$ $r < \max\{p,q\}$ if it is isotopic to a  $(p,q)$--torus knot
in which $r$ adjacent strands are twisted $n$ full twists. Consider the
twist knot $ K _m = T(7,17,2,5m - 3)$. The result then follows from the
following two theorems recalling an old result by Norwood, see ~\cite{No}.
He proves that a knot a with a fundamental group that is generated by
two elements is prime. Hence tunnel number one knots must be prime.

\begin{theorem}[Morimoto, Sakuma and Yokota ~\cite{MSY}]
For all integers $m$ the knots $K_m$ do not have a $(1,1)$--decomposition.
\end{theorem}

For tunnel number one knots we have the following theorem by Morimoto
~\cite{Mo3}:

\begin{theorem}[Morimoto ~\cite{Mo3}]
Let $K_1$ and $K_2$ be tunnel number one knots.  Then $K_1 \# K_2$
has tunnel number two if and only if at least one of $K_1$ and $K_2$
admits a $(1,1)$--decomposition.
\end{theorem}

More examples of  tunnel number one knots which do not have
$(1,1)$--decompositions have been discovered by Eudave--Mu\~noz in
~\cite{Eu}. They are quite complicated so we will not describe them
here. Recently Johnson and Thompson ~\cite{JT} have proved:

\begin{theorem}[Johnson and Thompson~\cite{JT}] \label{1nknots}
For every $n \in \naturals$, there is a tunnel  number one knot $K$
such that $K $ does not have a $(1, n)$ decomposition.
\end{theorem}

A stronger theorem which generalizes the above was also recently proved
by Minsky, Moriah and Schleimer ~\cite{MMS}:

\begin{theorem}[Minsky, Moriah and Schleimer ~\cite{MMS}] \label{gbknots} 
For any positive integers 
$t$ and $b$ there is a knot   $K \subset S^3$ with tunnel number $t$ so
that  $K$ has no $(t, b)$--decomposition. 
\end{theorem}

One should point out that both of the above theorems are not constructive.

The above theorem was proved utilizing the following theorem, which
is of independent interest, and work of Scharlemann--Tomova ~\cite{ST1}
and Tomova which generalizes her work in \cite{To}:

\begin{theorem}
For any pair of integers $g > 1$ and $n > 0$ there is a knot $K \subset
S^3$ and a genus $g$ splitting of $E(K)$ having distance greater than $n$.
\end{theorem}

Hence all of the above knots $K$ which are tunnel number one knots
must have $t(K \# K) = t(K) + t(K) + 1 = 3$.  The emergence of these
examples in recent years led to the following conjecture by
Morimoto~\cite{Mo2} and separately by the author ~\cite{Mr2}:

\begin{conjecture} 
\label{Conj:Morimoto}
If $K_1$ and $K_2$ are knots in $S^3$ then $t(K_1 \# K_2) = t(K_1) +
t(K_2) + 1$ if and only if neither of $K_1$ and $K_2$ are $\mu$--primitive.
\end{conjecture}

\begin{remark} \label{converse}
The conjecture is known to be true by work of Morimoto for tunnel number
one knots~\cite{Mo3} and for knots which are connected sums of two prime
knots each of which is also  $m$--small~\cite{Mo2}. It has been further
generalized by Kobayashi and Rieck  ~\cite{KR1} to general knots which
are $m$--small, that is, knots which do not contained essential surfaces
with meridional boundary components.
\end{remark}

However the situation is not that simple as there is some evidence for
the falsity of the conjecture:

Let  $nK$ denote  $K\#K\# \cdots \#K$, $n$ times.

\begin{enumerate}
\item It is a theorem of Kobayashi--Rieck ~\cite{KR2} that if there
exists a knot $K \subset S^3$ such that both $K$ and $2K = K \# K$
are not $\mu$--primitive then the above \fullref{Conj:Morimoto} is false.

\item It is also a theorem of theirs ~\cite{KR2} that if a knot
$K$ has a $(t,n)$--decomposition then the knot $nK$ is $\mu$--primitive.
\end{enumerate}

Note that all actual  known examples of knots which are not $(t,1)$
are $(t,2)$. Hence they fail as candidates for a counter example to the
above conjecture.

Note also that $(1)$ follows from $(2)$. Here is the proof given in
~\cite{KR2}: Let $m$ be the minimal positive integer  so that $t(mK)
< mt(K) + (m - 1)$.  Such an $m$ exists because of $(2)$ and $n \geq
2$. If $m = 2$ we are done as $K$ is the counter example.

Assume that $m \geq 3$. Set $K_1 = 2K$ and $K_2 = (m - 2) K$. By the
assumption and the minimality of $m$ we have $t(K_1) = 2t(K) + 1$ and
$t(K_2) =  (m - 2)t(K) + (m - 3)$.  However $t(K_1\#K_2) = t(mK) < mt(K)
+ (m - 1) = t(K_1) + t(K_2) + 1$. So $K_1$ and $K_2$ are a contradiction
to  \fullref{Conj:Morimoto}.

We give a proof for $(2)$ in \fullref{Add} below. Hence  they asked
the following question ~\cite[1.9]{KR2}: Are there knots $K \subset
S^3$ so that $t(K) = t$ and which admit a  $(t,n)$--decomposition with $n$
minimal  and $n \geq 3$?

\begin{remark}	
\label{Rem:Tricky}
Theorems~\ref{gbknots} and ~\ref{1nknots} give a positive answer to the
above question.  However showing that a knot $K$ which is a non-trivial
connected sum has a $(t(K),b)$--decomposition must be very tricky since
these knots contain essential tori so all their Heegaard splittings are
distance at most $2$ by Hempel~\cite{He}.  Hence the idea behind the proof of
Theorems~\ref{gbknots} and ~\ref{1nknots} will not work in this case
and one needs a completely different technique.
\end{remark}

Though the above \fullref{Rem:Tricky} shows that we are far away from
finding a counterexample to the Morimoto Conjecture using these techniques
it is  the belief of the author that the conjecture is false. In fact,
very recently the existence of counter examples has been announced by
Kobayashi and Rieck ~\cite{KR3,KR4}. A better conjecture
would be:

\begin{conjecture}
\label{Conj:??}  If $K_1, K_2$ are \emph{prime} knots in $S^3$ then
$t(K_1 \# K_2) = t(K_1) + t(K_2) + 1$ if and only if neither of $K_1$
and $K_2$ are $\mu$--primitive.
\end{conjecture}

\subsection{(2) Additive}\label{Add}

If $K_1 \subset S^3$ is $\mu$--primitive then for any other knot $K_2
\subset S^3$ we have $t(K) \leq t(K_1) + t(K_2)$.  What else can we say
about knots which are $\mu$--primitive?  The first thing to consider is
knots which have more than one primitive meridian.

\begin{definition}
We say that the knot $K \subset S^3$ is  {\it $n$--primitive} if  $E(K)$
has a minimal genus Heegaard splitting $(V,W)$ so that $\partial E(K)
\subset V$  and with $n$ primitive meridians corresponding to $n$
disjoint pairwise non-isotopic disks $D_1, \ldots,D_n \subset W$.
\end{definition}

\begin{example}
Let  $K_1$ be a knot in a $2n$--plat projection satisfying the following
conditions:

Let $B$ denote the $2n$--braid underlying the $2n$--plat projection
of $K$. Consider the image $\hat\rho(B)$ of $B$ under the Burau
representation of the $2n$--braid group, with the variable $t$ evaluated at
$-1$. Let $\alpha$ be the greatest common divisor of the entries  of the matrix $L \circ
\hat\rho(B) \circ M$. Assume  that $\alpha \neq 1$. Let $\beta$ be the
determinant of the matrix $N \circ \hat\rho(B) \circ M$.  Then there is
a  minimal genus $n$ Heegaard splitting of $E(K)$ with $n$ different
primitive meridians. In this case  we have of course that $n = t(K)
+ 1= g$. Furthermore if $K_2$ is another such knot then $t(K_1\# K_2)
= t(K_1) + t(K_2)$ (see Lustig--Moriah~\cite{LM3}). In other words the set of knots
which satisfy the conditions above is closed under connected sum.
\end{example} 

It is clear that if $K$ has $n$ primitive meridians then this procedure
can be iterated an arbitrary number of times with the corresponding
inequality. That is, if $K_1,  \ldots, K_n \subset S^3$ are any set of knots,
then $t(K \# K_1 \# \cdots \# K_n) \leq t(K) + t(K_1) + \cdots + t(K_n)$.

\begin{definition}\label{not-primitive}
If $ K \subset S^3 $ such that $t(K) = t$ does not have a primitive
meridian then it has a $ (t,n) $ decomposition for some minimal $
n $. We will say that $ K $ has a $ \frac{1}{n}$--\emph{primitive meridian}.
\end{definition}

In order to justify the above ``multiplicative" language we prove the
following theorem which is fact $(2)$ of \fullref{converse}.  This proof
is somewhat different than the one given in Kobayashi--Rieck~\cite{KR2}.

\begin{theorem}
Suppose the knot  $K$ has a $\frac{1}{n}$--primitive meridian; then $nK =
K\# \cdots \#K$ $n$--times  has a Heegaard splitting of genus $n t(K) +
n$  which has a $1$--primitive meridian. That is, $nK$ is $\mu$--primitive.
\end{theorem}

\begin{proof}
Let $K \subset S^3 $ be a knot with a $(t(K), n) $ decomposition. Then
there is a Heegaard splitting $(V_1, V_2)$ of $E(K)$ of genus $t(K)
+ n$ so that $\partial E(K) \subset V_1$.  It is obtained by taking
$n$ $1$--handles which are regular neighborhoods of the $n$ arcs $\{
t_1 ,\dots, t_n \} $ from one handlebody $W_2$ in the $(t(K), n )$
decomposition and adding them to the other handlebody $W_1$, then removing
a smaller regular neighborhood of $K$ from the modified $W_1$  to obtain a
compression body $V_1$ and a  handlebody  $V_2$ both of genus $t(K) + n$.

The boundary of the cocore disks of the tunnels now determines a
collection of $n$ simple closed curves $\{ \gamma_1 ,\dots, \gamma_n
\}  \subset \partial V_2$ which are by definition primitive.  They have
the additional property that the essential disks $\{ D_1 ,\dots, D_n \}
\subset V_2$ that each curve $\{ \gamma_1 ,\dots, \gamma_n \}  \subset
\partial V_2$ intersect in a single point can be chosen to be pairwise
disjoint and non-isotopic.

Now consider $(U^1_1,U^1_2), \dots , (U^{n-1}_1,U^{n-1}_2) $, $n-1$ copies
of a minimal genus Heegaard splitting  for  $E(K)$ so that   $\partial
E(K) \subset U_1 $. For each $i$ cut $U^i_1$  along a vertical annulus
$A^i$ so that $\partial A^i = \alpha_1 \cup \alpha_2$ with $\alpha_1
\subset \partial E(K)$.  This operation leaves two images $A_1^i, A_2^i$
of the annulus on the resulting handlebody. Similarly cut $V_1$ along
$n-1$ vertical annuli $B_i$ corresponding to the cocore disks of the
tunnels $\{ t_1 ,\dots, t_{n-1} \} $ . Attach a copy of the cut open
$U^i_1$ to the cut open $V_1$ by identifying the images  of $A_1^i $
and $B_1^i $ and  $A_2^i $ and $B_2^i$ to obtain a compression body.

Attach a copy of $U^i_2$ to $V_2$ by identifying the image of  an
annulus neighborhood of $\alpha_2$ in $\partial U^i_2 $ with the annulus
neighborhood of  $\gamma_i$. Since $\gamma_i$ is primitive we obtain
a handlebody of genus $  t(K) + n + (n-1)(t(K) +1) - (n-1)  = n t(K) +
n$  which determines  a Heegaard splitting $(H_1,H_2)$ of $E(nK)$. The
meridian corresponding to the $t_n$ arc is clearly primitive from the
construction meeting the essential disk $D_n$ of $H_2$ in a single point.
\end{proof} 

We now state some interesting  facts:

\begin{claim}
If a  knot $K \subset S^3$  is $\mu$--primitive then the corresponding
Heegaard splitting $(V,W)$ has distance $d(V,W) \leq 2$.
\end{claim}

\begin{proof}
We assume that $\partial E(K) \subset V$.  The compression body  $V$
has many disks which are  disjoint from the vertical annulus  $A$ given
by the assumption.  Denote one such disk by $D_0$. Let $D \subset W$
be the disk intersecting $A$ in a single point. There is an essential
disk $D' \subset W$ composed of two copies of $D$ and a band running
along the boundary component $a_1$ of $A$ contained in  $\partial _+  V$
which is disjoint from  $a_1$. Hence $d_{\mathcal{C}(S)}(D_0, a_1) = 1$
and $d_{\mathcal{C}(S)}(D', a_1) = 1$ so the Heegaard splitting $(V,W)$
has distance $d(V,W) \leq 2$.
\end{proof}  

In ~\cite{KR2} Kobayashi and Rieck make the following definition: 

\begin{definition}
The {\it growth rate of  tunnel number} is defined to be
$$t_{\text{gr}} (K) = \lim \frac {t(nK) - nt(K) } {n  - 1}.$$
It is also referred to as the {\it tunnel gradient}.
\end{definition}

\begin{proposition}
If $K$ is a knot on a $2n$--plat satisfying the conditions of Lustig and
Moriah ~\cite{LM3} then $t_{\text{gr}} (K) = 0$.
\end{proposition}

\begin{proof}
As $K$ satisfies the conditions of ~\cite{LM3} the tunnel number is
additive for connected sum of copies of $K$. It also has $n \geq 3$
primitive  meridians and by Proposition 6.2 of  ~\cite{LM3} the set of
knots which satisfy the  ~\cite{LM3} conditions is  closed under connected
sum. Hence we have $t(nK) = nt(K) $ so  $t_{\text{gr}} (K) =  \lim \frac {t(nK)
- nt(K) } {n  - 1} = 0$.
\end{proof}

\begin{remark}
Note that there are no known examples of  knots $K \subset S^3$ so that
$t(K \# K) < 2t(K)$. So we do not know of any knots with negative tunnel
gradient $t_{\text{gr}} (K)$.
\end{remark}

\subsection{(3) Sub-additive} \label{sub}

The first result showing that there are knot exteriors which behave in a
sub-additive way under connected sum is due to Morimoto in ~\cite{Mo1}. He
proved the following (also mentioned in \fullref{down}):

\begin{theorem}[Morimoto ~\cite{Mo1}]
Let   $K_1= K^n(-2,3,-3,2)$ be the twisted pretzel, where between the two
$3$--strands of twists one introduces  an odd number $n \in \mathbb{Z},
n \notin \{-1, 0, 1\}$ of ``horizontal" crossings, and $K_2$ is any
$2$--bridge knot.  Then  $t(K_1) = 2$  and $t(K_1 \# K_2) = 2 < 2 + 1 =
t(K_1) + t(K_2)$.
\end{theorem}

This result was then generalized by Kobayashi ~\cite{Ko5} who proved:

\begin{theorem}[Kobayashi ~\cite{Ko5}]
For each positive integer $m$, there exist knots $K_1, K_2$ so that
$t(K_1) + t(K_2) - t(K_1\#K_2) > m$.
\end{theorem}

Here is the idea of the proof:

The knot $K_1 = 2m K^n$ is  a connected sum of $2m$ copies of $K^n$,
the $n$ twisted pretzel $K^n(-2,3-3,2)$ as in \fullref{down}. The knot
$K_2 = 6m T(2,3)$ is  a connected sum of $6m$ copies of trefoils. It is
known by ~\cite{LM3} that $t(K_2) = 6m$ and Kobayashi proves, using the
tangle structure of  $K^n(-2,3-3,2)$, that $t(K_1) \geq 3m$. Now, using
the fact that connected sum is commutative, it follows from Morimoto's
~\cite{Mo1} ,above, that $t(K_1 \# K_2) \leq 8m$. Hence one obtains the
required result that $t(K_1) + t(K_2) - t(K_1\#K_2) > m$.

\begin{remark}
The ideas in the above proof are put together craftily to achieve large
degeneration of tunnel number; however, it is an iterated use of the
Morimoto example ~\cite{Mo1}.  There are no known other examples of two
knots which are sub-additive.
\end{remark}

There are some results which give some conditions and bounds for the
possible degeneration of tunnel number. The first \footnote{There is an
earlier slightly weaker result of Schultens ~\cite{Sch1} which states:
Theorem: For small knots $ K_1, K_2, t(K_1) + t(K_2) - 1 \leq t(K_1\#K_2)
\leq t(K_1)+t(K_2)+1 $.} such result was due to Morimoto and Schultens
~\cite{MS} who proved the following theorem:

Recall first that a knot $K \subset S^3$ is {\it small} if $E(K)$  does
not contain a closed essential (incompressible and not boundary parallel
) surface. It follows from Culler, Gordon, Luecke and Shalen
~\cite[Theorem~2.0.3]{CGLS} that if $E(K)$
contains an essential meridional (that is, a properly embedded bounded
surface so that all boundary components are meridians of $K$)  surface,
then $E(K)$ contains a closed essential surface.

\begin{theorem}[Morimoto--Schultens ~\cite{MS}]
If both of $K_1, K_2 \subset S^3$ are small then $t(K_1) +  t(K_2)
\leq t(K_1 \# K_2)$.
\end{theorem}

That is, the presence of meridional essential surfaces is required for
the degeneration of tunnel number.  This result was later improved  on
by Morimoto ~\cite{Mo4} who proved:

\begin{theorem}[Morimoto ~\cite{Mo4}]
Let $M_1, \dots, M_n$ be orientable closed $3$--manifolds which do not
have lens spaces summands. Let $K_1\subset M_1, \dots, K_n\subset M_n$
be  knots.  If none of  $E(K_i) = M_i - N(K_i)  , i = 1, \ldots, n$ contain
essential meridional surfaces then
$$t(K_1 \# \cdots \# K_n) \geq t(K_1) + \cdots + t(K_n).$$
\end{theorem}

In a slightly different flavor we have work of Scharlemann and Schultens:

\begin{theorem}[Scharlemann--Schultens~\cite{SS1}]
If $K _i \subset S^3, i = 1, \ldots, n$, are non-trivial  knots then $t(K_1
\# \cdots \# K_n) > n$.
\end{theorem}

A much stronger result is:

\begin{theorem}[Scharlemann--Schultens~\cite{SS2}]
Let $K_1, \ldots,  K_n$ be prime knots in $S^3$ then:
\begin{enumerate}
\item $t(K_1 \# \cdots \# K_n) \geq 1/3(t(K_1) + \cdots + t(K_n))$.
\item If none of the $K_i$ are $2$--bridge knots then
$$t(K_1 \# \cdots \# K_n)\geq 2/5(t(K_1) + \cdots + t(K_n).$$
\item $t(K_1 \#  K_2) \geq 2/5(t(K_1) + t(K_2))$.
\end{enumerate}
\end{theorem}

This leads to the following questions and conjecture:

\begin{remark}
We see from the above discussion that meridional essential surfaces are
required for the  tunnel number to degenerate. However it follows from
Lustig--Moriah ~\cite{LM3,LM4} that there are many knots which
contain many different  meridional essential surfaces but for which the
tunnel number is additive.
\end{remark}

Hence:

\begin{question}
What are the properties of meridional essential surfaces which ensure
that the tunnel number degenerate? Can these surfaces be classified?
\end{question}

\begin{question}
Are there knots which are not  $K_1 = K^n(-2,3,-3,2)$  and $2$--bridge
knots so that $t(K_1 \# K_2)  <   t(K_1) + t(K_2)$?
\end{question}

\begin{question}
Are there knots which are not composites of $K_1 = K^n(-2,3,-3,2)$
and $2$--bridge knots so that $t(K_1 \# K_2)  <   t(K_1) + t(K_2) - m$,
for $m \geq 1$?
\end{question}

\begin{conjecture}
If $K_1$ and  $K_2$ are {\it prime} knots for  which $t(K_1 \# K_2)
< t(K_1) + t(K_2) $ then $t(K_1 \# K_2)  =  t(K_1) +  t(K_2) - 1$.
\end{conjecture}

\bibliographystyle{gtart}
\bibliography{link}

\begin{thebibliography}{}
\providecommand\bibmarginpar{\leavevmode\marginpar}
\def\urlstyle#1{{\tt #1}}

\bibitem{Al}
\textbf{J Alexander}, \emph{A lemma on systems of knotted curves}, Proc. Nat.
  Acad. Sci. USA 9 (1963) 103--124

\bibitem{BCW}
\textbf{D Bachman}, \textbf{D Cooper}, \textbf{M\,E White},
  \href{http://dx.doi.org/10.2140/agt.2004.4.31} {\emph{Large embedded balls
  and {H}eegaard genus in negative curvature}}, Algebr. Geom. Topol. 4 (2004)
  31--47 \xox{MR}{2031911}

\bibitem{BM}
\textbf{S\,A Bleiler}, \textbf{Y Moriah},
  \href{http://dx.doi.org/10.1007/BF01456837} {\emph{Heegaard splittings and
  branched coverings of $B^3$}}, Math. Ann. 281 (1988) 531--543
  \xox{MR}{958258}

\bibitem{BO}
\textbf{M Boileau}, \textbf{J-P Otal}, \emph{Groupe des diff\'eotopies de
  certaines vari\'et\'es de {S}eifert}, C. R. Acad. Sci. Paris S\'er. I Math.
  303 (1986) 19--22 \xox{MR}{849619}

\bibitem{BRZ}
\textbf{M Boileau}, \textbf{M Rost}, \textbf{H Zieschang},
  \href{http://dx.doi.org/10.1007/BF01456287} {\emph{On {H}eegaard
  decompositions of torus knot exteriors and related {S}eifert fibre spaces}},
  Math. Ann. 279 (1988) 553--581 \xox{MR}{922434}

\bibitem{BZ}
\textbf{M Boileau}, \textbf{H Zieschang},
  \href{http://dx.doi.org/10.1007/BF01388469} {\emph{Heegaard genus of closed
  orientable {S}eifert {$3$}-manifolds}}, Invent. Math. 76 (1984) 455--468
  \xox{MR}{746538}

\bibitem{BZ1}
\textbf{G Burde}, \textbf{H Zieschang}, \emph{Knots}, de Gruyter Studies in
  Mathematics 5, Walter de Gruyter \& Co., Berlin (1985) \xox{MR}{808776}

\bibitem{CG}
\textbf{A\,J Casson}, \textbf{C\,M Gordon},
  \href{http://dx.doi.org/10.1016/0166-8641(87)90092-7} {\emph{Reducing
  {H}eegaard splittings}}, Topology Appl. 27 (1987) 275--283 \xox{MR}{918537}

\bibitem{CGLS}
\textbf{M Culler}, \textbf{C\,M Gordon}, \textbf{J Luecke}, \textbf{P\,B
  Shalen}, \href{http://dx.doi.org/10.2307/1971311} {\emph{Dehn surgery on
  knots}}, Ann. of Math. $(2)$ 125 (1987) 237--300 \xox{MR}{881270}

\bibitem{Do}
\textbf{H Doll}, \href{http://dx.doi.org/10.1007/BF01934349} {\emph{A
  generalized bridge number for links in $3$--manifolds}}, Math. Ann. 294
  (1992) 701--717 \xox{MR}{1190452}

\bibitem{Eu}
\textbf{M Eudave-Mu{\~n}oz}, \href{http://dx.doi.org/10.1142/S0218216506004804}
  {\emph{Incompressible surfaces and $(1,1)$--knots}}, J. Knot Theory
  Ramifications 15 (2006) 935--948 \xox{MR}{2251034}

\bibitem{Fo}
\textbf{R\,H Fox}, \emph{Free differential calculus I: Derivation in the free
  group ring}, Ann. of Math. 55 (1953) 547--560

\bibitem{Fu}
\textbf{K Funcke}, \href{http://dx.doi.org/10.1007/BF01247307} {\emph{Nicht
  frei \"aquivalente {D}arstellungen von {K}notengruppen mit einer
  definierierenden {R}elation}}, Math. Z. 141 (1975) 205--217 \xox{MR}{0372044}

\bibitem{GL}
\textbf{C\,M Gordon}, \textbf{J Luecke},
  \href{http://dx.doi.org/10.2140/agt.2006.6.2051} {\emph{Knots with unknotting
  number 1 and essential {C}onway spheres}}, Algebr. Geom. Topol. 6 (2006)
  2051--2116 \xox{MR}{2263059}

\bibitem{GR}
\textbf{C\,M Gordon}, \textbf{A\,W Reid},
  \href{http://dx.doi.org/10.1142/S0218216595000193} {\emph{Tangle
  decompositions of tunnel number one knots and links}}, J. Knot Theory
  Ramifications 4 (1995) 389--409 \xox{MR}{1347361}

\bibitem{Ha}
\textbf{Y Hagiwara}, \emph{Reidemeister--Singer distance for unknotting tunnels
  of a knot}, Kobe J. Math. 11 (1994) 89--100 \xox{MR}{1309994}

\bibitem{He}
\textbf{J Hempel}, \emph{3--manifolds}, Annals of Mathematics Studies 86,
  Princeton University Press (1976)

\bibitem{He1}
\textbf{J Hempel}, \href{http://dx.doi.org/10.1016/S0040-9383(00)00033-1}
  {\emph{3--manifolds as viewed from the curve complex}}, Topology 40 (2001)
  631--657 \xox{MR}{1838999}

\bibitem{Ja}
\textbf{W Jaco}, \emph{Lectures on three-manifold topology}, CBMS Regional
  Conference Series in Mathematics 43, American Mathematical Society,
  Providence, R.I. (1980) \xox{MR}{565450}

\bibitem{Jo}
\textbf{K Johannson}, \emph{Heegaard surfaces in {H}aken $3$--manifolds}, Bull.
  Amer. Math. Soc. $($N.S.$)$ 23 (1990) 91--98 \xox{MR}{1027902}

\bibitem{Jo1}
\textbf{K Johannson}, \emph{Topology and combinatorics of 3--manifolds},
  Lecture Notes in Mathematics 1599, Springer, Berlin (1995) \xox{MR}{1439249}

\bibitem{JT}
\textbf{J Johnson}, \textbf{A Thompson}, \emph{On tunnel number one knots which
  are not $(1,n)$}  \xox{arXiv}{math.GT/0606226}

\bibitem{KW}
\textbf{I Kapovich}, \textbf{R Weidmann}, \emph{Kleinian groups and the rank
  problem}, Geom. Topol.  (to appear)

\bibitem{Ko}
\textbf{T Kobayashi}, \emph{There exist 3-manifolds with arbitrarily high genus
  irreducible Heegaard splittings}, unpublished

\bibitem{Ko6}
\textbf{T Kobayashi}, \emph{Structures of the {H}aken manifolds with {H}eegaard
  splittings of genus two}, Osaka J. Math. 21 (1984) 437--455 \xox{MR}{752472}

\bibitem{Ko2}
\textbf{T Kobayashi}, \emph{A criterion for detecting inequivalent tunnels for
  a knot}, Math. Proc. Cambridge Philos. Soc. 107 (1990) 483--491
  \xox{MR}{1041480}

\bibitem{Ko1}
\textbf{T Kobayashi}, \emph{A construction of $3$--manifolds whose
  homeomorphism classes of {H}eegaard splittings have polynomial growth}, Osaka
  J. Math. 29 (1992) 653--674 \xox{MR}{1192734}

\bibitem{Ko5}
\textbf{T Kobayashi}, \href{http://dx.doi.org/10.1142/S0218216594000137}
  {\emph{A construction of arbitrarily high degeneration of tunnel numbers of
  knots under connected sum}}, J. Knot Theory Ramifications 3 (1994) 179--186
  \xox{MR}{1279920}

\bibitem{Ko4}
\textbf{T Kobayashi}, \href{http://dx.doi.org/10.2140/gt.2001.5.609}
  {\emph{Heegaard splittings of exteriors of two bridge knots}}, Geom. Topol. 5
  (2001) 609--650 \xox{MR}{1857522}

\bibitem{KR3}
\textbf{T Kobayashi}, \textbf{Y Rieck}, \emph{Knot exteriors with additive
  Heegaard genus and Morimoto's Conjecture}  \xox{arXiv}{math.GT/0701765}

\bibitem{KR4}
\textbf{T Kobayashi}, \textbf{Y Rieck}, \emph{Knots with $g(E(K))=2$ and
  $g(E(K\#K\#K))=6$ and Morimoto's Conjecture}  \xox{arXiv}{math.GT/070176}

\bibitem{KR}
\textbf{T Kobayashi}, \textbf{Y Rieck}, \emph{Manifolds admitting both strongly
  irreducible and weakly reducible Heegard splittings}, preprint

\bibitem{KR2}
\textbf{T Kobayashi}, \textbf{Y Rieck}, \emph{On the growth rate of tunnel
  number of knots}  \xox{arXiv}{math.GT/0402025}

\bibitem{KR1}
\textbf{T Kobayashi}, \textbf{Y Rieck}, \emph{Heegaard genus of the connected
  sum of $m$--small knots}, Comm. Anal. Geom.  (to appear)
  \xox{arXiv}{math.GT/0503229}

\bibitem{LM1}
\textbf{M Lustig}, \textbf{Y Moriah},
  \href{http://dx.doi.org/10.1016/0040-9383(91)90005-O} {\emph{Nielsen
  equivalence in {F}uchsian groups and {S}eifert fibered spaces}}, Topology 30
  (1991) 191--204 \xox{MR}{1098913}

\bibitem{LM3}
\textbf{M Lustig}, \textbf{Y Moriah},
  \href{http://dx.doi.org/10.1007/BF01444882} {\emph{Generalized Montesinos
  knots, tunnels and $\mathcal{N}$--torsion}}, Math. Ann. 295 (1993) 167--189
  \xox{MR}{1198847}

\bibitem{LM2}
\textbf{M Lustig}, \textbf{Y Moriah}, \emph{Generating systems of groups and
  Reidemeister--Whitehead torsion}, J. Algebra 157 (1993) 170--198

\bibitem{LM4}
\textbf{M Lustig}, \textbf{Y Moriah},
  \href{http://dx.doi.org/10.1007/PL00004347} {\emph{On the complexity of the
  {H}eegaard structure of hyperbolic {$3$}-manifolds}}, Math. Z. 226 (1997)
  349--358 \xox{MR}{1483536}

\bibitem{LM5}
\textbf{M Lustig}, \textbf{Y Moriah},
  \href{http://dx.doi.org/10.1016/S0166-8641(97)00232-0} {\emph{Closed
  incompressible surfaces in complements of wide knots and links}}, Topology
  Appl. 92 (1999) 1--13 \xox{MR}{1670164}

\bibitem{LM6}
\textbf{M Lustig}, \textbf{Y Moriah},
  \href{http://dx.doi.org/10.1016/S0040-9383(99)00024-5} {\emph{3--manifolds
  with irreducible Heegaard splittings of high genus}}, Topology 39 (2000)
  589--618 \xox{MR}{1746911}

\bibitem{LM7}
\textbf{M Lustig}, \textbf{Y Moriah},
  \href{http://dx.doi.org/10.1016/j.top.2004.01.004} {\emph{A finiteness result
  for {H}eegaard splittings}}, Topology 43 (2004) 1165--1182 \xox{MR}{2079999}

\bibitem{LS}
\textbf{R\,C Lyndon}, \textbf{P\,E Schupp}, \emph{Combinatorial group theory},
  Ergebnisse der Mathematik und ihrer Grenzgebiete 89, Springer, Berlin (1977)
  \xox{MR}{0577064}

\bibitem{MMS}
\textbf{Y Minsky}, \textbf{Y Moriah}, \textbf{S Schleimer}, \emph{Knots of high
  distance}, preprint

\bibitem{M}
\textbf{Y Moriah}, \emph{Heegaard splittings and group presentations}, PhD
  thesis, University of Texas at Austin (1986)

\bibitem{Mr0}
\textbf{Y Moriah}, \href{http://dx.doi.org/10.1007/BF01388781} {\emph{Heegaard
  splittings of {S}eifert fibered spaces}}, Invent. Math. 91 (1988) 465--481
  \xox{MR}{928492}

\bibitem{Mr1}
\textbf{Y Moriah}, \href{http://dx.doi.org/10.1016/S0166-8641(01)00303-0}
  {\emph{On boundary primitive manifolds and a theorem of Casson--Gordon}},
  Topology Appl. 125 (2002) 571--579 \xox{MR}{1935173}

\bibitem{Mr2}
\textbf{Y Moriah}, \emph{Connected sums of knots and weakly reducible Heegaard
  splittings}, Topology Appl. 141 (2004) 1--20 \xox{MR}{2058679}

\bibitem{MR}
\textbf{Y Moriah}, \textbf{H Rubinstein}, \emph{Heegaard structures of
  negatively curved $3$--manifolds}, Comm. Anal. Geom. 5 (1997) 375--412
  \xox{MR}{1487722}

\bibitem{MSS}
\textbf{Y Moriah}, \textbf{S Schleimer}, \textbf{E Sedgwick},
  \href{http://projecteuclid.org/getRecord?id=euclid.cag/1175790071}
  {\emph{Heegaard splittings of the form $H{+}nK$}}, Comm. Anal. Geom. 14
  (2006) 215--247 \xox{MR}{2255010}

\bibitem{MS1}
\textbf{Y Moriah}, \textbf{J Schultens},
  \href{http://dx.doi.org/10.1016/S0040-9383(97)00072-4} {\emph{Irreducible
  {H}eegaard splittings of {S}eifert fibered spaces are either vertical or
  horizontal}}, Topology 37 (1998) 1089--1112 \xox{MR}{1650355}

\bibitem{MS2}
\textbf{Y Moriah}, \textbf{E Sedgwick},
  \href{http://dx.doi.org/10.1142/S0218216504003470} {\emph{Closed essential
  surfaces and weakly reducible {H}eegaard splittings in manifolds with
  boundary}}, J. Knot Theory Ramifications 13 (2004) 829--843 \xox{MR}{2088748}

\bibitem{Mo3}
\textbf{K Morimoto}, \href{http://dx.doi.org/10.1016/0166-8641(93)90099-Y}
  {\emph{On the additivity of tunnel number of knots}}, Topology Appl. 53
  (1993) 37--66 \xox{MR}{1243869}

\bibitem{Mo}
\textbf{K Morimoto}, \href{http://dx.doi.org/10.1016/0166-8641(94)90099-X}
  {\emph{On composite tunnel number one links}}, Topology Appl. 59 (1994)
  59--71 \xox{MR}{1293117}

\bibitem{Mo1}
\textbf{K Morimoto}, \href{http://dx.doi.org/10.2307/2161103} {\emph{There are
  knots whose tunnel numbers go down under connected sum}}, Proc. Amer. Math.
  Soc. 123 (1995) 3527--3532 \xox{MR}{1317043}

\bibitem{Mo2}
\textbf{K Morimoto}, \href{http://dx.doi.org/10.1007/PL00004411} {\emph{On the
  super additivity of tunnel number of knots}}, Math. Ann. 317 (2000) 489--508
  \xox{MR}{1776114}

\bibitem{Mo4}
\textbf{K Morimoto}, \href{http://dx.doi.org/10.1016/S0040-9383(98)00070-6}
  {\emph{Tunnel number, connected sum and meridional essential surfaces}},
  Topology 39 (2000) 469--485 \xox{MR}{1746903}

\bibitem{MS3}
\textbf{K Morimoto}, \textbf{M Sakuma},
  \href{http://dx.doi.org/10.1007/BF01446565} {\emph{On unknotting tunnels for
  knots}}, Math. Ann. 289 (1991) 143--167 \xox{MR}{1087243}

\bibitem{MSY}
\textbf{K Morimoto}, \textbf{M Sakuma}, \textbf{Y Yokota}, \emph{Examples of
  tunnel number one knots which have the property ``$1+1=3$''}, Math. Proc.
  Cambridge Philos. Soc. 119 (1996) 113--118 \xox{MR}{1356163}

\bibitem{MS}
\textbf{K Morimoto}, \textbf{J Schultens},
  \href{http://dx.doi.org/10.1090/S0002-9939-99-05160-6} {\emph{Tunnel numbers
  of small knots do not go down under connected sum}}, Proc. Amer. Math. Soc.
  128 (2000) 269--278 \xox{MR}{1641065}

\bibitem{Ni}
\textbf{J Nielsen}, \href{http://dx.doi.org/10.1007/BF01457113} {\emph{Die
  {I}somorphismen der allgemeinen, unendlichen {G}ruppe mit zwei
  {E}rzeugenden}}, Math. Ann. 78 (1964) 385--397 \xox{MR}{1511907}
  \xox{JFM}{46.0175.01}

\bibitem{No}
\textbf{F\,H Norwood}, \href{http://dx.doi.org/10.2307/2044414} {\emph{Every
  two-generator knot is prime}}, Proc. Amer. Math. Soc. 86 (1982) 143--147
  \xox{MR}{663884}

\bibitem{OS}
\textbf{P Ozsv{\'a}th}, \textbf{Z Szab{\'o}},
  \href{http://dx.doi.org/10.1016/j.top.2005.01.002} {\emph{Knots with
  unknotting number one and {H}eegaard {F}loer homology}}, Topology 44 (2005)
  705--745 \xox{MR}{2136532}

\bibitem{Ro}
\textbf{D Rolfsen}, \emph{Knots and links}, Mathematics Lecture Series 7,
  Publish or Perish, Berkeley, CA (1976) \xox{MR}{0515288}

\bibitem{RS}
\textbf{H Rubinstein}, \textbf{M Scharlemann},
  \href{http://dx.doi.org/10.1016/0040-9383(95)00055-0} {\emph{Comparing
  {H}eegaard splittings of non-{H}aken $3$--manifolds}}, Topology 35 (1996)
  1005--1026 \xox{MR}{1404921}

\bibitem{Sa}
\textbf{M Sakuma}, \emph{Manifolds with infinitely many non-isotopic Heegaard
  splittings of minimal genus}, from: ``(Unofficial) Proceedings of the
  Conference on Various Structures on Knots and their Applications (Osaka City
  University)'' (1988)  172--179

\bibitem{Sc}
\textbf{M Scharlemann}, \href{http://dx.doi.org/10.1090/S0002-9947-03-03182-9}
  {\emph{There are no unexpected tunnel number one knots of genus one}}, Trans.
  Amer. Math. Soc. 356 (2004) 1385--1442 \xox{MR}{2034312}

\bibitem{SS1}
\textbf{M Scharlemann}, \textbf{J Schultens},
  \href{http://dx.doi.org/10.1016/S0040-9383(98)00002-0} {\emph{The tunnel
  number of the sum of {$n$} knots is at least {$n$}}}, Topology 38 (1999)
  265--270 \xox{MR}{1660345}

\bibitem{SS2}
\textbf{M Scharlemann}, \textbf{J Schultens},
  \href{http://dx.doi.org/10.1090/S0002-9947-00-02654-4} {\emph{Comparing
  {H}eegaard and {JSJ} structures of orientable 3--manifolds}}, Trans. Amer.
  Math. Soc. 353 (2001) 557--584 \xox{MR}{1804508}

\bibitem{ST}
\textbf{M Scharlemann}, \textbf{A Thompson}, \emph{Heegaard splittings of
  $(\text{surface}){\times}I$ are standard}, Math. Ann. 295 (1993) 549--564

\bibitem{ST1}
\textbf{M Scharlemann}, \textbf{M Tomova},
  \href{http://dx.doi.org/10.2140/gt.2006.10.593} {\emph{Alternate {H}eegaard
  genus bounds distance}}, Geom. Topol. 10 (2006) 593--617 \xox{MR}{2224466}

\bibitem{Sh}
\textbf{H Schubert}, \emph{Knoten mit zwei Br\"ucken}, Math. Z. 65 (1956)
  133--170

\bibitem{Sch}
\textbf{J Schultens}, \emph{The classification of {H}eegaard splittings for
  (compact orientable surface){$\,\times\, S\sp 1$}}, Proc. London Math. Soc.
  $(3)$ 67 (1993) 425--448 \xox{MR}{1226608}

\bibitem{Sch1}
\textbf{J Schultens}, \href{http://dx.doi.org/10.1007/s000140050131}
  {\emph{Additivity of tunnel number for small knots}}, Comment. Math. Helv. 75
  (2000) 353--367 \xox{MR}{1793793}

\bibitem{SW}
\textbf{J Schultens}, \textbf{R Weidmann}, \emph{On the geometric and algebraic
  rank of graph manifolds}, preprint

\bibitem{Se}
\textbf{E Sedgwick}, \href{http://dx.doi.org/10.1007/s002080050064} {\emph{An
  infinite collection of {H}eegaard splittings that are equivalent after one
  stabilization}}, Math. Ann. 308 (1997) 65--72 \xox{MR}{1446199}

\bibitem{Se2}
\textbf{E Sedgwick}, \emph{The irreducibility of {H}eegaard splittings of
  {S}eifert fibered spaces}, Pacific J. Math. 190 (1999) 173--199
  \xox{MR}{1722770}

\bibitem{Se1}
\textbf{E Sedgwick}, \href{http://dx.doi.org/10.2140/agt.2001.1.763}
  {\emph{Genus two 3--manifolds are built from handle number one pieces}},
  Algebr. Geom. Topol. 1 (2001) 763--790

\bibitem{So}
\textbf{J Souto}, \href{http://dx.doi.org/10.1016/j.top.2004.09.001} {\emph{A
  note on the tameness of hyperbolic 3--manifolds}}, Topology 44 (2005)
  459--474 \xox{MR}{2114957}

\bibitem{To}
\textbf{M Tomova}, \emph{Multiple bridge surfaces restrict knot distance}
  \xox{arXiv}{math.GT/0511139}

\end{thebibliography}

\end{document}